\documentstyle[12pt]{article}
\newlength{\abstractwidth}
\renewcommand{\thefootnote}{\fnsymbol{footnote}}
\renewcommand{\thanks}[1]{\footnote{#1}} 
\newcommand{\starttext}{ \setcounter{footnote}{0}
\renewcommand{\thefootnote}{\arabic{footnote}}}

\newcommand{\be}{\begin{equation}}
\newcommand{\bea}{\begin{eqnarray}}
\newcommand{\eea}{\end{eqnarray}} \newcommand{\ee}{\end{equation}}
 \newcommand{\<}{\langle}
\renewcommand{\>}{\rangle} \def\ba{\begin{eqnarray}}
\def\ea{\end{eqnarray}}

\def\o{\omega}

\def\tr{{\rm tr}}
\def\det{{\rm det}}

\def\log{\,{\rm log}\,}

\def\o{\omega}

\def\e{\varepsilon}

\def\o{\omega}

\def\na{\nabla}

\def\ge{\geq}
\def\le{\leq}

\def\p{\partial}

\def\na{{\nabla}}

\def\[{{\bf [}}
\def\]{{\bf ]}}

\def\ddbar{i\p\bar\p}

\def\mathbb{\bf}

\def\eqref{\ref}

\newcommand{\neweqref}[1]{(\ref{#1})}

\begin{document}
\starttext \baselineskip=18pt \setcounter{footnote}{0}
\newtheorem{theorem}{Theorem}
\newtheorem{lemma}{Lemma}
\newtheorem{corollary}{Corollary}
\newtheorem{definition}{Definition}
\newtheorem{conjecture}{Conjecture}
\newtheorem{proposition}{Proposition}

\begin{center}
{\Large \bf AUXILIARY MONGE-AMPERE EQUATIONS IN GEOMETRIC ANALYSIS
\footnote{Contribution to Surveys in Differential Geometry in honor of S.S. Chern.
Work supported in part by the National Science Foundation under grant DMS-22-03273.}}

\medskip
\centerline{Bin Guo and  Duong H. Phong}

\medskip

\begin{abstract}

{\footnotesize This is an introduction to a particular class of auxiliary complex Monge-Amp\`ere equations which had been instrumental in $L^\infty$ estimates for fully non-linear equations and various questions in complex geometry. The essential comparison inequalities are reviewed and shown to apply in many contexts. Adapted to symplectic geometry, with the auxiliary equation given now by a real Monge-Amp\`ere equation, the method gives an improvement of an earlier theorem of Tosatti-Weinkove-Yau, reducing Donaldson's conjecture on the Calabi-Yau equation with a taming symplectic form from an exponential bound to an $L^1$ bound.}

\end{abstract}

\end{center}

\baselineskip=15pt
\setcounter{equation}{0}
\setcounter{footnote}{0}

\section{Introduction}
\setcounter{equation}{0}

It is well-known that comparisons with an auxiliary equation can be of powerful assistance in the study of partial differential equations. Comparisons with harmonic functions were used early on by De Giorgi \cite{DeG}, and have been applied since in many contexts, including e.g. in the book of L. Simon \cite{Sl}, and in X.J. Wang's method for Schauder estimates for the Poisson equation \cite{W}. Since then, comparisons with many other equations have proved to be effective. Complex Monge-Amp\`ere equations are particularly suitable as auxiliary equations, since the existence and smoothness of their solutions for given right hand sides have been established by Yau \cite{Y} in the case of compact K\"ahler manifolds, and by Caffarelli, Kohn, Nirenberg, and Spruck \cite{CKNS} in the case of the Dirichlet problem. Notable successes of Monge-Amp\`ere equations as auxiliary equation include the estimates of Dinew and Kolodziej \cite{DK1} relating volume and capacity, the bounds of Song and Tian \cite{ST} for the K\"ahler-Ricci flow, and the entropy estimates of Chen and Cheng \cite{CC} for the constant scalar curvature equation.

\medskip
Very recently, a specific class of auxiliary Monge-Amp\`ere equations has been instrumental in many significant advances in complex geometry.
This class was first introduced 
by the authors in joint work with F. Tong in \cite{GPT}, and already led to a pure PDE proof of the $L^\infty$ estimates of Kolodziej \cite{K}, a goal which had eluded researchers in the field for close to a quarter of a century. As a PDE proof, the method extends immediately Kolodziej's estimates to a general class of non-linear equations satisfying a structural condition. Remarkably, this class has been shown by Harvey and Lawson \cite{HL3} to be quite large, and include in particular all invariant Garding-Dirichlet operators. But the method turns out to be even more flexible and powerful than naively anticipated, and it has had since many unexpected applications. These include stability estimates for Monge-Amp\`ere and Hessian equations
\cite{GPT1}; $L^\infty$ estimates for Monge-Amp\`ere equations on nef classes rather than just K\"ahler classes \cite{GPTW1}; sharp modulus of continuity for non-H\"older solutions \cite{GPTW2}; extensions to parabolic equations \cite{CC1}; extensions to equations on Hermitian manifolds \cite{GP}; extensions to form-type equations \cite{GP}; lower bounds for the Green's function \cite{GPS}; uniform entropy estimates \cite{GPa}; and diameter estimates and convergence theorems in K\"ahler geometry not requiring bounds on the Ricci curvature \cite{GS,GPSS}. 

\medskip

The main purpose of this paper is to provide a survey of these developments. We shall describe in some detail the essential features of the particular class of auxiliary complex Monge-Amp\`ere equations of interest. The main applications are then sketched, with explanations of how the auxiliary Monge-Amp\`ere equations are used and precise statements of the results obtained. In all but one case, the full treatment is left to references to the original papers in the literature. The one exceptional case is the application to taming symplectic forms on almost-K\"ahler manifolds. It is a conjecture of Donaldson \cite{D}, motivated by symplectic geometry, that on a compact $4$-manifold equipped with an almost-complex structure $J$ and a taming symplectic form $\Omega$, the Calabi-Yau equation would admit a priori bounds to all orders. 
This conjecture had been reduced by Weinkove \cite{W} to an $L^\infty$ bound for the potential, and subsequently by Tosatti, Weinkove, and Yau \cite{TWY} to a single exponential estimate.
Using an auxiliary real Monge-Amp\`ere equation, we can reduce it further to a single $L^1$ estimate. This result is treated in detail because it is new and does not appear anywhere else. But it may also be noteworthy as evidence that the methods here can extend to the real or symplectic context as well.

\section{The auxiliary Monge-Amp\`ere equation}
\setcounter{equation}{0}

We begin by describing the class of Monge-Amp\`ere equations which will serve later as comparison equations. Some key inequalities needed for the maximum principle are broadly described, which can be adapted later for different applications.

\medskip

Let $(X,\o_X)$ be a compact $n$-dimensional Hermitian manifold. For any Hermitian form $\o$ on $X$ and any smooth function $\varphi$ with ${\rm sup}_X\varphi=0$, set $\o_\varphi=\o+i\p\bar\p\varphi$, and consider the relative endomorphism 
\bea
h_\varphi=\o_X^{-1}\o_\varphi.
\eea
We denote by $\lambda[h_\varphi]$ the (un)-ordered vector of its eigenvalues.
Let $f(\lambda)$ be a function on a convex cone $\Gamma\subset {\bf R}^n$ invariant under permutations, 
and consider the family of equations parametrized by $\o$,
\bea
\label{eq:f}
f(\lambda[h_\varphi])=k_\o(z),\quad \lambda[h_\varphi]\in\Gamma,
\eea
where $k_\o(z)$ is a given positive function. For some estimates, it is convenient to introduce the constant $c_\o>0$ defined by the following normalization 
\bea
\label{rho}
k(z)=c_\o\,e^{F_\o},
\quad
\int_X e^{nF_\o}\o_X^n=\int_X\o_X^n.
\eea

As in \cite{GPT}, we require that $f:\Gamma\to {\bf R}_+$ satisfies

 (1) $\Gamma\subset {\mathbb R}^n$ is a symmetric cone with 
\begin{equation}\label{eqn:cone}
\Gamma_n\subset \Gamma \subset \Gamma_1.
\end{equation}
Here $\Gamma_k$ is the cone of vectors $\lambda$ with $\sigma_j(\lambda)>0$ for $1\leq j\leq k$, where $\sigma_j(\lambda)$ is the $j$-th symmetric polynomial in $\lambda$. In particular, $\Gamma_1$ is the half-space defined by $\lambda_1+\cdots+\lambda_n>0$, and $\Gamma_n$ is the first octant, defined by $\lambda_j>0$ for $1\leq j\leq n$.

(2) $f(\lambda)$ is symmetric in $\lambda = (\lambda_1,\ldots, \lambda_n)\in \Gamma$ and it is homogeneous of degree one;

 (3) $\frac{\partial f}{\partial \lambda_j}>0$ for each $j=1,\ldots, n$ and $\lambda\in \Gamma$;

 (4) There is a $\gamma>0$ such that 
\begin{equation}\label{eqn:structure}
\prod_{j=1}^n \frac{\partial f}{\partial \lambda_j}\ge \gamma\,{\rm det}^{-1}(\o_X)\,\quad \forall \lambda\in \Gamma.
\end{equation}

It is well-known that equations such as the Monge-Amp\`ere equation, with $f(\lambda)=(\prod_{j=1}^n\lambda_j)^{1\over n}$, or the Hessian equation with $f(\lambda)=\sigma_k(\lambda)^{1\over k}$ where $\sigma_k$ is the $k$-th order symmetric polynomial, or the $p$-Monge-Amp\`ere equation of Harvey and Lawson \cite{HL1,HL2} with
$$f(\lambda) = \Big( \prod_{I} \lambda_I\Big)^{\frac{n!}{(n-p)!p!}}$$
 where $I$ runs over all distinct multi-indices $1\le {i_1}<\cdots < {i_p}\le n$, $\lambda_I = \lambda_{i_1} + \cdots + \lambda_{i_p}$, and $\Gamma$ is the cone defined by $\lambda_I>0$ for all  $p$-indices $I$, all satisfy the structural condition (4). A remarkable recent result of Harvey and Lawson \cite{HL3} is that the condition (4) actually holds for very large classes of non-linear operators, including all invariant Garding-Dirichlet operators. As noted in \cite{HL3}, the condition (4) also arose independently in \cite{AO} in the study of $W^{2,p}$ interior regularity.

\smallskip
We would like to compare the solution $\varphi$ of the equation (\ref{eq:f}) with the solution $\psi$ of the following complex Monge-Amp\`ere equation
\bea
\label{eq:MA}
(\o+i\p\bar\p\psi)^n = {\tau (-\varphi+q(z)-s)\over A}k_\o^n(z)\,\o_X^n.
\qquad
{\rm sup}_X\psi=0
\eea
where $s\ge  0$ is a nonnegative constant, $\tau(t)$ is a smooth strictly positive function, $q(z)$ is a given function, and $A$ is a normalizing constant.

\medskip
Let $G^{j\bar k}$ be the linearized operator of $\log f(\lambda)$, defined by
\bea
G^{j\bar k}={\p\over \p h_{\bar kj}}\log f(\lambda[h])
\eea
Fix a point $z_0\in X$. Since $\o_X$ is positive definite, and $h_\varphi$ is a self-adjoint endomorphism with respect to $\o_X$,
we can choose a holomorphic coordinate system centered at $z_0$ where $(\o_X)_{\bar jk}=\delta_{jk}$ and  $h_\varphi(z_0)$ is diagonal, $(h_\varphi)_{\bar jk}(z_0)=\lambda_j\delta_{jk}$.
 In particular $(\o_\varphi)_{\bar kj}(z_0)=\lambda_j\delta_{jk}$.

\begin{lemma}
\label{lm:trace}
The linearized endomorphism $G^{j\bar k}$ then satisfies the following: for $\lambda\in \Gamma$

{\rm (a)} $G^{j\bar k}(z_0)={1\over f(\lambda)}{\p f\over\p \lambda_j}\delta_{jk}$, and $G^{j\bar k}$ is positive definite;

{\rm (b)} $G^{j\bar k}(\o_\varphi)_{\bar kj}(z_0)=1$ and $G^{ j\bar k}\o_{\bar kj}(z_0)\geq 0$, 

{\rm (c)} and, taking into account the equations satisfied by $\varphi$ and $\psi$,
\bea
{1\over n}G^{j\bar k}(\o_\psi)_{\bar kj}
\geq
({\gamma \tau(-\varphi+q-s)\over A})^{1\over n}.
\eea
\end{lemma}

\medskip
\noindent
{\it Proof.} The formula in (a) is a classic formula for the linearization of a fully non-linear operator $\log f(\lambda[h])$ at a diagonal matrix $h$ from the theory of non-linear equations \cite{S}. 
The positive-definiteness of $G^{j\bar k}$ is an immediate consequence of the ellipticity condition $\p_jf(\lambda)>0$ of the function $f$. Next, we can write at $z_0$,
\bea
G^{j\bar j}(\o_\varphi)_{\bar kj}=\sum_{j=1}^n {1\over f(\lambda)}{\p f\over \p\lambda_j}\lambda_j=1
\eea
by Euler's relation for homogeneous functions $f(\lambda)$ of degree $1$. This proves the first equation in (b). To establish the second equation in (b), we observe that since $(\o_X)_{\bar kj}=\delta_{jk}$ at $z_0$, $\o$ can be identified with the relative endomorphism $h_0=(\o_X)^{-1}\o$, which is hermitian and positive with respect to $\o_X$. The expression $G^{j\bar k}\o_{\bar kj}$ can then be identified with the trace ${\rm Tr}(Gh_0)$, which is then the inner product of two positive matrices. As such, it is positive, as can be seen for example by writing it in a basis where both endomorphisms are diagonal.

Finally, to establish (c), we apply the arithmetic-geometric inequality to write
\bea
{1\over n}G^{j\bar k}(\o_\psi)_{\bar kj}
&\geq& \bigg({\rm det}\,(G^{j\bar k}(\o_\psi)_{\bar km})\bigg)^{1\over n}
=\bigg({\rm det} G^{j\bar k}\cdot{\rm det}(\o_\psi)_{\bar km})\bigg)^{1\over n}
\nonumber\\
&=&\bigg({1\over f(\lambda)^n}(\prod_{j=1}^n{\p f\over \p\lambda_j})\cdot {\tau (-\varphi+q-s)\over A }k_\o(z)^n{\rm det}(\o_X)\bigg)^{1\over n}.
\eea
Applying the equation for $\varphi$, this inequality simplifies to
\bea
{1\over n}G^{j\bar k}(\o_\psi)_{\bar kj}
\geq 
\big(\prod_{j=1}^n{\p f\over \p\lambda_j}\cdot {\tau (-\varphi+q-s)\over A }{\rm det}(\o_X)\big)^{1\over n}
\eea
Applying now the structural condition on $f(\lambda)$, we obtain the desired inequality (c). Q.E.D.

\medskip
Next, we need a ``comparison function $\Phi$", relating the solution $\varphi$ of the equation (\ref{eq:f}) to the solution $\psi$ of the auxiliary Monge-Amp\`ere equation (\ref{eq:MA}). For the applications considered in this survey, the comparison function $\Phi$ is usually of the form
\bea
\label{Phi}
\Phi=-\e(-\psi+q(z)+\Lambda)^b-\varphi+\tilde q(z)-s
\eea
with $0<b<1$ a fixed constant, $\e,\Lambda$ non-negative constants to be chosen later, and $\tilde q(z)$ a smooth function.
We obtain the key inequality relating $\varphi$ and $\psi$ if we can choose all the data in $\Phi$ so as to guarantee that $\Phi\leq 0$ everywhere. To do so, we typically apply the maximum principle, and make use of the following general calculations:

\begin{lemma}
\label{lm:maximum}
Fix a point $z_0$ and a holomorphic coordinate system centered at $z_0$ as above. Then we have the following inequality at $z_0$
\bea
\label{eq:maximum}
G^{j\bar k}\Phi_{\bar kj}
&\geq&
\e nb(-\psi+q(z)+\Lambda)^{b-1}({\gamma \tau(-\varphi+q-s)\over A})^{1\over n}-1
\nonumber\\
&&
+G^{j\bar k}\big\{\o_{\bar kj}+(1-\e b(-\psi+q+\Lambda)^{b-1})(\o_q)_{\bar kj}+\tilde q_{\bar kj}\big\}.
\eea

\end{lemma}

\medskip
\noindent
{\it Proof.} A direct calculation gives
\bea
\Phi_{\bar kj}&=&\e b(-\psi+q(z)+\Lambda)^{b-1}(\psi-q)_{\bar kj}-\varphi_{\bar kj}+\tilde q_{\bar kj}\nonumber\\
&&-\e b(b-1)(-\psi+q(z)+\Lambda)^{b-2}(\psi-q)_j(\psi-q)_{\bar k}
\eea
Rewriting this expression using $(\psi-q)_{\bar kj}=(\o_\psi)_{\bar kj}-(\o_q)_{\bar kj}$ and $\varphi_{\bar kj}=(\o_\varphi)_{\bar kj}-\o_{\bar kj}$, we obtain
\bea
G^{j\bar k}\Phi_{\bar kj}
&=&\e b(-\psi+q+\Lambda)^{b-1}G^{j\bar k}(\o_\psi)_{\bar kj}
-\e b(-\psi+q+\Lambda)^{b-1}G^{j\bar k}(\o_q)_{\bar kj}
\nonumber\\
&&
-G^{j\bar k}(\o_\varphi)_{\bar kj}+G^{j\bar k}\o_{\bar kj}+G^{j\bar k}\tilde q_{\bar kj}
\nonumber\\
&&-\e b(b-1)(-\psi+q+\Lambda)^{b-2}G^{j\bar k}(\psi-q)_j(\psi-q)_{\bar k}
\eea
Since $G^{j\bar k}>0$ and $0<b<1$, we can drop from the right hand side the expression involving $G^{j\bar k}(\psi-q)_j(\psi-q)_{\bar k}$. The desired inequality follows then from the inequalities in Lemma \ref{lm:trace}. Q.E.D.

\bigskip
The lemma is particularly useful when the last term on the right hand side of (\ref{eq:maximum}) happens to be positive and can be dropped. 
We obtain then an upper bound for the solution $\varphi$ of the given equation (\ref{eq:f}) in terms of the solution $\psi$ of the auxiliary Monge-Amp\`ere equation.
The following is the simplest illustration, which will be shown later to apply to $L^\infty$ estimates on K\"ahler manifolds:

\begin{lemma}
\label{lm:use}
Let $\varphi$ and $\psi$ satisfy the equations (\ref{eq:f}) and (\ref{eq:MA}), under the preceding hypotheses on the operator $f(\lambda)$. Assume that the function $\tau(t)$ satisfies the condition
\bea
\label{tau}
\tau(t)\geq t^a,
\quad t\in [0,\infty)
\eea
for some fixed power $a>0$. Then for any $s\geq 0$, we have
\bea
-\varphi-s\leq c_{n,a,\gamma}A^{1\over a+n}(-\psi+\Lambda)^{n\over n+a}
\eea
for all $z\in X$ and all $s\geq 0$, if the constants $b$, $\e$ and $\Lambda$ are chosen to be 
\bea
b={n\over n+a},
\quad
\e=(nb\gamma^{1\over n})^{-{n\over a+n}}A^{1\over a+n},
\quad
\Lambda^{1-b}= \e b.
\eea
Here $c_{n,a,\gamma}$ is a constant depending only on $n,a,\gamma$.
\end{lemma}

\medskip
\noindent
{\it Proof.} Let $\Phi$ be defined as in (\ref{Phi}), with $q(z)=\tilde q(z)=0$. We apply Lemma \ref{lm:maximum}. Since $q=\tilde q=0$, we have $\o_q=\o$ and the last expression on the right hand side of (\ref{eq:maximum}) reduces to
\bea
G^{j\bar k}\big\{\o_{\bar kj}+(1-\e b(-\psi+q+\Lambda^{b-1}))(\o_q)_{\bar kj}+\tilde q_{\bar kj}\big\}
&=& G^{j\bar k}\big\{1-\e b(-\psi+\Lambda)^{b-1}\big\}\o_{\bar kj}
\nonumber\\
&\geq & (1-\e b\Lambda^{b-1})G^{j\bar k}\o_{\bar kj}\geq 0\label{eqn:2.19old}
\eea
since $-\psi\geq 0$, $0<b<1$, and $G^{j\bar k}\o_{\bar kj}\geq 0$ by (b) of Lemma \ref{lm:trace}.

Let $z_0$ be a point where $\Phi$ attains its maximum on $X$. We shall show that $\Phi(z_0)\leq 0$. 
If $-\varphi(z_0)-s\leq 0$, the function $\Phi$ is manifestly $\leq 0$ at its maximum $z_0$, and we are done. 
Otherwise, we note that $0\geq G^{j\bar k}\Phi_{\bar kj}(z_0)$ and apply Lemma \ref{lm:maximum}. As just noted, 
we can drop the last term on the right hand side of (\ref{eq:maximum}), and bound $\tau(-\varphi-s)$ from below by $(-\varphi-s)^a$. We find
\bea
1\geq n\e b(-\psi+\Lambda)^{b-1}({\gamma (-\varphi-s)^a\over A})^{1\over n}
\eea
at $z_0$, which can be rewritten as
\bea
-(nb\gamma^{1\over n})^{-{n\over a+n}}A^{1\over a+n}(-\psi+\Lambda)^{n\over n+a}-\varphi-s\leq 0.
\eea
With the choice of $b$, $\e$, and $\Lambda$ indicated in the lemma, we can recognize the left hand side as $\Phi(z_0)$. Since $z_0$ is a maximum for $\Phi$, it follows that $\Phi(z)\leq 0$ for any $z\in X$.
This last statement can be recast in the form stated in the Lemma. Q.E.D.

\section{Application to the compact K\"ahler case}
\label{section 3}
\setcounter{equation}{0}

The first application that we discuss is the one where the above class of auxiliary equations was originally introduced \cite{GPT}, in order to provide a PDE proof of the $L^\infty$ estimates for the complex Monge-Amp\`ere equation originally established by Kolodziej \cite{K}. 
We discuss this case in some detail, since it also serves as a template for other subsequent applications.

\medskip

Assume in this section that $X$ is a compact K\"ahler manifold, and both $\o_X$ and $\o$ are K\"ahler forms. We would like to derive $L^\infty$ estimates for the solutions $\varphi$ of the equation (\ref{eq:f}). 
Our goal is to find in the present case of non-linear equations an adaptation of the strategy going back to De Giorgi for $L^\infty$ bounds for linear equations in divergence form. In this strategy, a lower bound for $\varphi$ is obtained by showing that the set
\bea
\Omega_s=\{\varphi<-s\}
\eea
is empty starting from some $S_0$ which can be estimated. We shall deduce this from suitable growth conditions on the function
\bea\label{eqn:phi}
\phi(s)={1\over V_\omega }\int_{\Omega_s}k_\o^n(z)\o_X^n, 
\quad
V_\o=\int_X\o^n,
\eea
which will follow themselves from a reverse H\"older inequality for the key function $A_s$ defined for $s>0$ by
\bea
\label{def:As}
A_s={1\over V_\o}\int_{\Omega_s}(-\varphi-s)k_\o^n(z)\o_X^n.
\eea

For this we need an auxiliary Monge-Amp\`ere equation. Let $\tau_\ell(t)$ a sequence of smooth strictly positive functions on ${\bf R}$ which decreases monotonically to the function ${\bf R}\ni t\to t\chi_{{\bf R}_+}(t)$ as $\ell\to\infty$, and which are uniformly bounded from above by $1+t\chi_{{\bf R}_+}(t)$. Here $\chi_{{\bf R}_+}$ is the characteristic function of ${\bf R}_+$.
For each $s\in {\bf R}$, let $\psi_\ell(z)$ be the solution of the following auxiliary Monge-Amp\`ere equation
\bea
(\o+i\p\bar\p\psi_\ell)^n
=
{\tau_\ell(-\varphi-s)\over A_{\ell,s}}k_\o^n(z)
\o_X^n, 
\eea
where the constant $A_{\ell,s}$ is defined  by
\bea
A_{\ell,s}={1\over V_\o}\int_X\tau_\ell(-\varphi-s)k_\o^n(z)
\o_X^n.
\eea
By Yau's theorem \cite{Y}, the above equation admits a unique smooth and $\o$-plurisubharmonic (PSH) solution $\psi_\ell$ normalized by ${\rm sup}_X\psi_\ell=0$.

\medskip
We can now apply Lemma \ref{lm:use}, with $a=1$, $b=n/(n+1)$, $q\equiv 0$. Since $\psi_\ell\leq 0$, the condition $-\psi_\ell+q\geq 0$ is satisfied, and we obtain an estimate of the form
\bea
{-\varphi-s\over A_{\ell,s}^{1/ (n+1)}}
\leq c_{n,\gamma}(-\psi_\ell+c_{n,\gamma}A_{\ell,s})^{n\over n+1}
\eea
where $c_{n,\gamma}$ denote generically positive constants depending only on $n$ and $\gamma$. Restricting to the set $\Omega_s=\{-\varphi-s>0\}$, we can take the $n/(n+1)$-root of both sides, multiply the resulting inequality by a constant $\beta_0$, taking the exponential, and integrate over $\Omega_s$. We obtain
\bea
\int_{\Omega_s}{\rm exp}\big\{\beta_0({-\varphi-s\over A_{\ell,s}^{1/( n+1)}})^{n+1\over n}\big\}\o_X^n
\leq {\rm exp}(c_{n,\gamma}\beta_0A_{\ell,s})\int_{\Omega_s}{\rm exp}(-c_{n,\gamma}\beta_0\psi_\ell)\o_X^n.
\eea
We can now invoke the well-known inequality for $\alpha$-invariants, which states that for all $\o$ K\"ahler forms with $\o\leq \kappa \o_X$ for some fixed $\kappa>0$, there exists a constant $\alpha$ so that,
\bea
\int_X{\rm exp}(-\alpha_0\psi)\o_X^n
\leq C(\alpha_0,n,\o_X,\kappa)
\eea
for any $\alpha_0<\alpha$ and any $\o$-plurisubharmonic function $\psi$ with ${\rm sup}_X\psi=0$ \cite{H, T}. Thus we have
\bea
\int_{\Omega_s}{\rm exp}\big\{\beta_0({-\varphi-s\over A_{\ell,s}^{1/( n+1)}})^{n+1\over n}\big\}\o_X^n
\leq c_{\alpha_0,n,\o_X,\kappa,\gamma}{\rm exp}({c_{n,\gamma,\beta_0}A_{\ell,s}}).
\eea
We can now let $\ell\to\infty$ and obtain
\bea
\label{As}
\int_{\Omega_s}{\rm exp}\big\{\beta_0({-\varphi-s\over A_{s}^{1/( n+1)}})^{n+1\over n}\big\}\o_X^n
\leq c_{\alpha_0,n,\o_X,\kappa,\gamma}{\rm exp}({c_{n,\gamma,\beta_0}A_{s}}).
\eea
where $A_s={\rm lim}_{\ell\to\infty}A_{\ell,s}$ with $A_s$ defined as in (\ref{def:As}).

\smallskip
Let $E$ be the energy, defined by 
\bea
\label{energy}
E={1\over V_\o}\int_X (-\varphi)k_\o^n(z)\o_X^n.
\eea 
Then $A_s\leq E$ for any $s\geq 0$, and the preceding inequality implies the following,
\bea
\label{As1}
\int_{\Omega_s}{\rm exp}\big\{\beta_0({-\varphi-s\over A_{s}^{1/( n+1)}})^{n+1\over n}\big\}\o_X^n
\leq c\,{\rm exp}(c\,E)
\eea
which suffices for our purposes.

\medskip

The inequality (\ref{As}) is the key inequality in our method. It is not difficult to show that it implies a reverse H\"older inequality. For this, we apply the Young's inequality in the following form
\bea
UV\leq U\eta(U)+V\eta^{-1}(V)
\eea
for any monotone strictly increasing function $\eta:{\bf R}_+\to{\bf R}_+$ with ${\rm lim}_{u\to 0}\eta(u)=0$. Here $\eta^{-1}$ is the inverse of the function $\eta$. We make the choice $\eta(u)=(\log(1+u))^p$, $\eta^{-1}(v)={\rm exp}(V^{1\over p})-1$, and $U=e^{nF_\o}$, $V=v(z)^p$, where we have rewritten $k_\o(z)$ as $k_\o(z)=c_\o e^{F_\o}$ as in (\ref{rho}). This gives
\bea
e^{nF_\o}v(z)^p
\leq e^{nF_\o}\log^p(1+e^{nF_\o})+v(z)^p(e^{v(z)}-1)
\leq
c_p\big\{e^{nF_\o}(1+|nF_\o|^p)+e^{2v(z)}).
\nonumber
\eea
Next, take 
\bea
v(z)={1\over 2}\beta_0\Big({-\varphi-s\over A_s^{1/(n+1)}}\Big)^{n+1\over n}
\eea
and integrate both sides over $\Omega_s$. In view of (\ref{As}), we find

\begin{lemma}
\label{reverseHolder}
The following inequality holds
\bea
{1\over V_\o}\int_{\Omega_s}(-\varphi-s)^{p{n+1\over n}}k_\o^n(z)\o_X^n
\leq 
c{c_\o^n\over V_\o}\,A_s^{p\over n}
(\|e^{nF_\o}\|_{L^1(\log L)^p}+ce^{cE})
\eea
where $c=C(X,\o_X,\beta_0,p)$ is a constant.
\end{lemma}

We observe that $A_s$ is essentially the $L^1$ norm of $-\varphi-s$, so the inequality we just obtained can be interpreted as a reverse H\"older inequality. It implies immediately the following growth rate for the function $\phi(s)$ as in (\ref{eqn:phi}),

\begin{lemma}
\label{growth}
Fix $p>n$. Then we have the inequalities

{\rm (a)} $A_s\leq B_0\,\phi(s)^{1+\delta_0}$, for $\delta_0={(p-n)/np}>0$.

{\rm (b)} For any $r>0$, $A_s\geq r\phi(r+s)$,

{\rm (c)} $r\phi(s+r)\leq B_0\,\phi(s)^{1+\delta_0}$,\\
where the constant $B_0= [c{c_\omega^n\over V_\omega}(\|e^{nF_\o}\|_{L^1(\log L)^p}+ce^{cE})]^{1\over p}$.

\end{lemma}

\noindent
{\it Proof.} We begin with the proof of (a). By H\"older's inequality, we have
\bea
A_s\leq \big\{{1\over V_\o}\int_{\Omega_s}(-\varphi-s)^{p{n+1\over n}}k_\o^n(z)\o_X^n\big\}^{n\over (n+1)p}({1\over V_\o}\int_{\Omega_s}k_\o^n(z)\o_X^n)^{1\over q}
\eea
where $q$ is defined by ${n\over p(n+1)}+{1\over q}=1$. The first factor on the right hand side has been shown to be bounded by a multiple of $A_s^{1/(n+1)}$. Putting this factor on the left hand side yields
\bea
A_s\leq B_0 ({1\over V_\o}\int_{\Omega_s}e^{nF_\o}\o_X^n)^{1+n\over qn}=B_0\phi(s)^{1+n\over qn}.
\eea
 The exponent $1+n/(qn)$ is readily worked out to be $1+\delta_0$, establishing (a). To see (b), we observe that, trivially, $\Omega_{r+s}=\{\varphi<-s-r\}\subset \Omega_s$, and hence
\bea
A_s={1\over V_\o}\int_{\Omega_s}(-\varphi-s)e^{nF_\o}\o_X^n
\geq 
r{1\over V_\o}\int_{\Omega_{s+r}}e^{nF_\o}\o_X^n=r\phi(s+r)
\eea
establishing (b). The last statement (c) is a trivial consequence of (a) and (b). Q.E.D.

\medskip
We can now invoke a classic lemma of De Giorgi, which says that positive monotone decreasing functions $\phi(s)$ which tend to $0$ as $s\to\infty$, and satisfy the growth rate condition stated as (c) in Lemma \ref{growth}, 
must vanish for some $s\geq S_0$, where $S_0>0$ can be estimated in terms of $B_0$ and $\delta_0$. But this implies that $\Omega_s$ must be empty for $s\geq S_0$, and hence $\varphi\geq -S_0$. Thus we have obtained the following bound for $\varphi$ \cite{GPT}:

\begin{theorem}
\label{General}
Let $\varphi $ a solution to the equation where the operator $f(\lambda)$ satisfies the structural conditions listed in (1-4). Assume that the K\"ahler form $\o$ satisfies the condition $\o\leq \kappa\,\o_X$ for some positive constant $\kappa$. Fix $p>n$. Then we have the $L^\infty$ bound
\bea
\varphi\geq -C
\eea
where $C$ is a constant depending only on $X,\o_X,n,p,\gamma,\kappa$ and the following three quantities
\bea
{c_\o^n\over V_\o}, \quad E, \quad {\rm Ent}_p(\o) =\|e^{nF_\o}\|_{L^1(\log L)^p}.
\eea
\end{theorem}

In the case of the Monge-Amp\`ere equation $f(\lambda)=(\prod_{j=1})^n\lambda_j)^{1\over n}$, it is easy to see that ${c_\o^n\over V_\o}={1\over [\o_X^n]}$ and that the energy $E$ can be bounded by a constant depending only on $X,\o_X,n,\gamma,\kappa$ and ${\rm Ent}_{\{p=1\}}(\o)$. Thus 
Theorem \ref{General} gives $L^\infty$ estimates for $\varphi$ depending only on the entropy ${\rm Ent}_p(\o)$ for $p>n$, recovering in this way the $L^\infty$ estimates of Kolodziej \cite{K}, even in the more general version allowing degenerations of the background metric established by Demailly-Pali \cite{DP}, and Eyssidieux-Guedj-Zeriahi \cite{EGZ}. We stress however that it holds for the general class of operators $f(\lambda)$ satisfying the structural conditions (1-4), which is quite large, as shown by Harvey and Lawson \cite{HL3}.

\section{Application to energy estimates from entropy}
\setcounter{equation}{0}

The previous Theorem \ref{General} had reduced $L^\infty$ bounds for general non-linear equations of the form (\ref{eq:f}) to the three quantities $c_\o^n/[\o^n]$, $E$, and ${\rm Ent}_p(\o)$ for $p>n$. Actually, it has been known for some time that, for fixed background metric $\o$, bounds for the energy $E$ can be derived from bounds for the entropy \cite{B, CC}. This was even one of the key steps in the work of X.X. Chen and J.R. Cheng \cite{CC} on the equation for K\"ahler metrics of constant scalar curvature. As shown in \cite{GPT}, the arguments of \cite{CC} can be adapted to provide bounds for the energy, assuming that the entropy is bounded. However, all of these bounds depend on the background metric $\o$, and are only useful for fixed $\o$. 

\medskip
This problem of bounding uniformly the energy by the entropy is addressed by the following theorem \cite{GP}:

\begin{theorem}\label{thm:main1}\label{thm:main}
Let $(X,\o_X)$ is a compact $n$-dimensional K\"ahler manifold without boundary. Let $\o$ be any K\"ahler form on $X$ with
\bea
\label{kappa}
\o\leq \kappa\,\o_X
\eea
for some constant $\kappa>0$. Consider the equation \neweqref{eq:f} with the operator $f(\lambda)$ satisfying the conditions (1-4).
Then for any $p>0$, any $C^2$ solution $\varphi$ to \neweqref{eq:f} satisfies the following 

{\rm (i)} Trudinger-like inequalities 
\begin{equation}\label{eqn:inequality}
\int_X e^{\alpha (-\varphi)^{q}} \omega_X^n\le C_T,
\end{equation} 

{\rm (ii)} and energy-like estimates
\begin{equation}\label{eqn:inequality 2}
\int_X (-\varphi)^{pq} e^{nF_\o} \omega_X^n\le C_e.
\end{equation}

Here the exponent $q$ is given by $q = \frac{n}{n-p}$ if $p<n$, and can be any strictly positive exponent if $p \geq n$. The constants $C_T$ and $C_e$ are computable constants depending only on $n, p, q, \omega_X, \kappa,\gamma$, and upper bounds for the following two quantities
\bea
{c_\o^n\over V_\omega},
\quad
{\mathrm{Ent}}_p({\o}) = \int_X e^{nF_\o} |F_\o|^p\omega_X^n,
\eea
and the term $\alpha>0$ is a constant that depends only on $n, p, \gamma, {c_\o^n\over V_\omega}$ and $\kappa$.

\end{theorem}

\medskip
We observe that, in the case of a fixed background K\"ahler metric $\o$, this theorem was proved as Theorem 3 in \cite{GPT}. As stressed above, the point of the new theorem is to have uniform estimates, even as the background metric $\o$ may degenerate to the boundary of the K\"ahler cone. For this same reason, we consider the case $p>n$. When $p>n$, when the background K\"ahler form $\o$ is fixed, it follows from \cite{GPT}, Theorem 1, that the solution $\varphi$ of the equation is actually bounded, and the above Trudinger-like and energy-like estimates follow at once. But here again, the existing results do not give the estimates uniform in $\o$ that we seek.

\medskip
We only describe the version of the auxiliary complex Monge-Amp\`ere equation that we need, together with the comparison inequality, leaving fuller details to \cite{GPa}.

\medskip
Thus let $a = pq = \frac{np}{n-p}$ ($a$ is {\em any} positive number if $p=n$). We solve the following complex Monge-Amp\`ere equation
\begin{equation}\label{eqn:aux}
(\omega + \ddbar \psi_{s,\ell})^n = \frac{\tau_\ell(-\varphi - s) ^a}{A_{s,\ell}} c_\o^n e^{nF_\o} \omega_X^n,\quad {\rm sup}_X \psi_{s,\ell} = 0.
\end{equation}
Here the constant $A_{s,\ell}$ is defined by 
\begin{equation}\label{eqn:Ask}
A_{s,\ell} = \frac{c^n}{V_\omega} \int_X \tau_\ell(-\varphi - s) ^ a e^{nF_\o} \omega_X^n
\end{equation}
to make the equation \neweqref{eqn:aux} compatible.  By assumption $[\omega]$ is a K\"ahler class, so by Yau's theorem \cite{Y}, \neweqref{eqn:aux} admits a unique smooth solution $\psi_{s,\ell}$. We observe that as $\ell\to\infty$
\begin{equation}\label{eqn:As}
A_{s,\ell}\to A_s: = \frac{c_\o^n}{V_\omega} \int_{\Omega_s} (-\varphi - s)^a e^{nF_\o} \omega_X^n. 
\end{equation}

\smallskip

Next we consider comparison functions of the form
$$\Phi := -\varepsilon (-\psi_{s,\ell} + \Lambda)^b - \varphi - s, $$
where $\varepsilon$ and $\Lambda$ are constants, and $b$ is a suitable power. Using the identities in Section 2 and the maximum principle, we then show that
\bea
\Phi\leq 0
\eea
if the constants are given by \begin{equation}\label{eqn:constant 1}
b = \frac{n}{n+a}\in(0,1), \quad \mbox{and }\varepsilon = \frac{1}{\gamma^{1/(n+a)} (n b)^{n/(n+a)}} A_{s,\ell}^{\frac{1}{n+a}},
\end{equation}
and $\Lambda$ is chosen so that $\varepsilon b \Lambda^{-(1-b)} = 1$, that is,
\begin{equation}\label{eqn:constant 2}
\Lambda = \frac{b^{1/(1-b)}}{(\gamma^{1/(n+a)} (nb)^{n/(n+a)})^{1/(1-b)}} A_{s,\ell}^{\frac{1}{a}}.
\end{equation}
The inequality $\Phi\le 0$ on $X$ implies that $-\varphi-s$ can be controlled by $c_{a,n,\gamma}A_{s,\ell}^{1\over n+a}
(-\psi_{s,\ell}+A_{s,\ell}^{1\over a})^b$, after which the theorem can be established following the template of Section 3.

\bigskip

\section{Application to stability estimates}
\setcounter{equation}{0}

It is not difficult to adapt the same method of $L^\infty$ estimates to stability estimates.
As in Section \ref{section 3}, let $(X,\omega_X)$ be a compact K\"ahler manifold, and $\omega$ be a K\"ahler metric such that $\omega \le \kappa \omega_X$ for some $\kappa>0$. We consider the following complex Monge-Amp\`ere equations
\begin{equation}\label{eqn:sbMA}
(\omega+ \ddbar u)^n = c_\omega e^f \omega_X^n, \mbox{ and } (\omega + \ddbar v)^n = c_\omega e^h \omega_X^n,
\end{equation}
with the constant $c_\omega = \int_X \omega^n$ and $\int_X e^f \omega_X^n = \int_X e^h \omega_X^n = 1$. The functions $u,v$ are normalized such that 
$${\rm max}_X (u - v) = {\rm max}_X (v - u).$$
Assume that for some $K>0$ and $p>n$,
$$\max(\| e^f\|_{L^1(\log L)^p (X,\omega_X)}, \| e^h\|_{L^1(\log L)^p (X,\omega_X)})\le K. $$
\begin{theorem}
\label{thm:GPTa}
Under these assumptions, there is a constant $C>0$ depending on $n,p, \omega_X, \kappa$, and $ K$ such that 
$${\rm sup}_X |u - v| \le C \| e^f - e^h\|_{L^1}^\beta,$$
where $\beta = \beta(n,p) = (n+3 + \frac{p-n}{pn})^{-1}>0$.
\end{theorem}

\smallskip

We outline the main idea of the proof, and refer to \cite{GPT1} for the details. Fix a small $r>0$. First we may assume $\int_{\{v\le u\}} (e^f + e^h)\omega_X^n \le 1$. We then consider the auxiliary equation
\bea
(\o +i\p\bar\p\psi)^n
= c_\omega {\tau_\ell(-\varphi +q(z)-s)\over A_{\ell,s}} e^h \o_X^n,\quad {\rm sup}_X \psi = 0,
\eea
with $\varphi=v-(1-r)u$, $q(z)=-3\beta_0r$, where $\beta_0$ is an upper bound of $\| u\|_{L^\infty}$ and $\| v\|_{L^\infty}$. $\beta_0$ exists from Section \ref{section 3}. Take the comparison function
\bea
\Phi=-\e(-\psi+\Lambda)^{n\over n+1}-\varphi+q(z)-s
\eea
with $\varepsilon = c_1(n) A_{\ell,s}^{1/(n+1)}$ and $\Lambda = c_2(n ) (\varepsilon/r) ^{n+1}$ for appropriate constants $c_1(n)$ and $c_2(n)$. Arguing as in Section \ref{section 3}, we arrive at $ v- u \ge - C r$ for a constant $C>0$ depending only on $n, p,\kappa, \omega_X$ and $K>0$. 

\smallskip

We remark that analogous stability estimates as in Theorem \ref{thm:GPTa} hold for complex Hessian equations as well
\cite{GPT1,DK2}.

\section{Application to the nef case}
\setcounter{equation}{0}

So far, we have considered equations whose background form $\o$ is K\"ahler (although it may degenerate to the boundary of the K\"ahler cone). In this section, following \cite{GPTW1}, we show how it can be adapted to nef classes.

\medskip

As above, $(X,\omega_X)$ is a K\"ahler manifold. We assume only the class $[\omega]$ is K\"ahler and $\omega\le \kappa \omega_X$, but the $(1,1)$-form $\omega$ may not be positive. 

We consider both the Monge-Amp\`ere and Hessian equations. While $L^\infty$ estimates do not hold in the usual form, the following estimates can be established \cite{GPTW1}:

\begin{theorem}
\label{nef}
Assume that $(X,\o_X)$ is an $n$-dimensional compact K\"ahler manifold, and let $\o$ be a closed $(1,1)$-form which may not be positive, but $[\o]$ is a K\"ahler class with $\o\leq\kappa\,\o_X$. 

{\rm (a)} Let $f(\lambda)=(\prod_{j=1}^n\lambda_j)^{1\over n}$, corresponding to the Monge-Amp\`ere equation. 
Recall that the {\em envelope} associated to $\omega$ is the $\o$-PSH function defined by
\bea
{\mathcal V}_\omega = \sup\{u|~ u\in PSH(X,\omega),\, u\le 0\}.
\eea
Then for any fixed $p>n$, there is a constant $C$ depending only on $\o_X,n,p,\kappa, {\rm Ent}_p(\o)$, so that the solution $\varphi$ of the equation (\ref{eq:f}) satisfies
\bea
0\leq -\varphi+{\mathcal V}_\o\leq C.
\eea

{\rm (b)} Let $f(\lambda)=\sigma_k(\lambda)^{1\over k}$, corresponding to the $k$-th Hessian equation. Define the envelope ${\mathcal V}_{\o,k}$ corresponding to the $\Gamma_k$ cone by
\bea
{\mathcal V}_{\o,k}={\rm sup} \,\{v\leq 0\}
\eea
where $v$ runs over non-negative $C^2$ functions with the vector of eigenvalues of the relative endomorphism $\o_X^{-1}(\o+i\p\bar\p v)$ lying in the $\Gamma_k$ cone. Define the energy $E_k$ by
\bea
E_k=\int_X (-\varphi+{\mathcal V}_{\o,k})e^{kF_\o}\o_X^n.
\eea
Recall also the constant $c_\o$ as defined for the equation (\ref{eq:f}). 
Then for any $p>n$, there is a constant $C$ depending only on $\o_X,n,o,\kappa$, $E_k$, and $\|e^{{ k}F_\o}\|_{L^1(\log L)^p}$, ${c_\o^n\over [\o]^k[\o_X]^{n-k}}$  so that
\bea
0\leq -\varphi+{\mathcal V}_{\o,k}\leq C.
\eea
\end{theorem}

We note that Part (a) on the Monge-Amp\`ere equation had been proved with different methods by Boucksom et al \cite{BEGZ} and Fu-Guo-Song \cite{FGS}. Part (b) on Hessian equations is new and due to \cite{GPTW1}.

\smallskip
As in earlier applications, we indicate only the auxiliary Monge-Amp\`ere equation to be used, and the comparison function $\Phi$. We only discuss the proof of Part (a), the one of Part (b) being similar.

\smallskip
A first technical difficulty is that ${\mathcal V}_\omega$ is only $C^{1,1}$. 
However, it has been shown by Berman \cite{B} that there exists a sequence of smooth and strictly $\omega$-PSH functions $\{u_\beta\}_{\beta=1}^\infty$ which converge uniformly to ${\mathcal V}_\omega$. In the following we can use $u_\beta$ in place of ${\mathcal V}_\omega$, and then take limit $\beta\to \infty$.

\medskip

The auxiliary Monge-Amp\`ere equation is then
\bea
(\o+i\p\bar\p\psi)^n={\tau_\ell(-\varphi+u_\beta-s)\over A_{\ell,\beta}}k_\o^n(z)\o_X^n,\quad {\rm sup}_X\psi = -1
\eea
and the comparison function $\Phi$ is defined by
\bea
\Phi=-\e(-\psi+u_\beta+1+\Lambda)^{n\over n+1}-\varphi+u_\beta-s.
\eea
In the general setting of Section 2, this corresponds to our general ansatz with $q=\tilde q = u_\beta$ and the simple shift in the normalization of $\psi$, with ${\rm sup}_X\psi=-1$. We note that an analogue of (\ref{eqn:2.19old}) still holds by the fact that $\omega_{u_\beta} >0$ although $\omega$ may not be positive.
We can then show that $\Phi\leq 0$, and establish the desired bounds following our general template.

\section{Application to the modulus of continuity}
\setcounter{equation}{0}

It has been shown by Kolodziej \cite{K} that the solution of the complex Monge-Amp\`ere equation on a compact K\"ahler manifold is H\"older continuous if the right hand side is of class $L^q$ for some $q>1$. He has also shown that if the solution is in some Orlicz space, then it must be continuous. But even so, his arguments are not direct, and we don't have any information on the modulus of continuity of the solution. The modulus of continuity is typically a delicate question. However, it turns out that it can also be addressed using the class of auxiliary Monge-Amp\`ere equations discussed in the present paper.

\medskip
More specifically, 
let $(X,\omega_X)$ be again a compact K\"ahler manifold of complex dimension $n$. We consider the complex Monge-Amp\`ere equation with $\int_X e^F \omega_X^n = \int_X \omega_X^n$
\begin{equation}\label{eqn:MAm}
(\omega_X + \ddbar \varphi)^n = e^F \omega_X^n, \quad \omega_\varphi = \omega_X + \ddbar \varphi>0. 
\end{equation}
Then the following estimate is established in \cite{GPTW2}:

\begin{theorem}
\label{thm:GPTW2}
Fix $p>n$. Then we have 
$$|\varphi (x) - \varphi(y)|\le \frac{C}{|\log d(x,y)|^\alpha},\quad \forall x,y\in X$$
for some constant $C>0$ depending on $n, p, \omega_X$ and $\| e^F\|_{L^1(\log L)^p}$. Here $d(x,y)$ denotes the geodesic distance of $x,y$ in the Riemannian manifold $(X,\omega_X)$, and $\alpha =\min \{\frac{p}{n+1}, \frac{p-n}{n}\}>0$.
\end{theorem}

The proof of Theorem \ref{thm:GPTW2} relies on the following auxiliary complex Monge-Amp\`ere equation:
\bea
(\o_X+i\p\bar\p\psi)^n={\tau_\ell(-\varphi+q(z)-s)\over A_{\ell,s}}e^F\o_X^n
\eea
with $q = (1- |\log \delta|^{-p/(n+1)}) \varphi_\delta - 2\delta $ and $\varphi_\delta$ is the (rescaled) Kiselman-Legendre transform of $\varphi$ at level $\delta>0$. With the comparison function
\bea
\Phi=-\e(-\psi+\Lambda)^{n\over n+1}-\varphi+q(z)-s,
\eea
we can argue as in Section \ref{section 3} to conclude that $\varphi_ \delta - \varphi \le 2\delta + |\log \delta|^{-p/(n+1)} \varphi_\delta + S_\infty$ for $S_\infty = {C}{|\log \delta|^{-(p-n)/n}}$. The proof of Theorem \ref{thm:GPTW2} then follows from the fact that $\varphi_\delta(z)$ is equal to $\max \varphi$ over the ball $B(z, \delta)$ up to a controlled error term. 

\section{Application to the Hermitian case}
\setcounter{equation}{0}

A striking feature of the above method of auxiliary Monge-Amp\`ere equations is the ease with which it can be adapted to the case of Hermitian manifolds. Thus we let $(X,\o_X)$ be a compact Hermitian manifold with $\o_X$ a fixed Hermitian metric, $\o_\varphi=\o_X+i\p\bar\p\varphi$, 
$h_\varphi=\o^{-1}_X\o_\varphi$, and consider the equation
\bea
\label{eq:f hermitian}
f(\lambda[h_\varphi])=e^{F},
\qquad \lambda[h_\varphi]\in \Gamma,
\qquad {\rm sup}_X\varphi=0,
\eea 
where the operator $f(\lambda)$ defined on a cone $\Gamma$ satisfies the structural conditions (1-4) spelled out in section \S 2. We have then \cite{GP}

\begin{theorem}
For any $p>n$, and $C^2$ solution $\varphi$ of the equation (\ref{eq:f hermitian}) must satisfy the $L^\infty$ bound
\bea
{\rm sup}_X|\varphi|\leq C
\eea
where $C$ is a constant depending only on $X,\o_X,n,p,\gamma$ and 
$
\|e^{nF}\|_{L^1(\log L)^p}$.
\end{theorem}

In the case of the Monge-Amp\`ere equation, $L^\infty$ estimates were first obtained in the Hermitian setting by Cherrier \cite{Ch} and Tosatti and Weinkove \cite{ToWe}, assuming a pointwise bound for $e^F$. The sharper version with entropy bounds $\|e^{nF}\|_{L^1(\log L)^p}$ was obtained by Dinew and Kolodziej \cite{DK1}, and required a highly non-trivial extension of pluripotential theory to the Hermitian setting. An approach based on envelopes has been recently proposed by Guedj and Lu \cite{GL}. Our theorem applies to much more general classes of equations, and its proof is arguably the simplest, as it bypasses the complicated integration by parts with torsion terms which arise in Hermitian geometry. We provide a brief sketch.

\medskip
A first observation is that the previous auxiliary complex Monge-Amp\`ere equation cannot be applied as it is. The reason is that, on a compact Hermitian manifold, unlike on K\"ahler manifolds, solutions exist only up to an undetermined constant. Because of this, we shall use instead as auxiliary equation the Dirichlet problem for a complex Monge-Amp\`ere equation on a Euclidean ball, which has been shown by Caffarelli, Kohn, Nirenberg, and Spruck \cite{CKNS} to admit always a smooth solution.

\medskip
Thus, fix $r_0$ small enough, but depending only on $(X,\o_X)$ so that, for any $z_0\in X$, there is a coordinate system $z$ centered at $z_0$ so that
\bea
{1\over 2}i\p\bar\p |z|^2\leq\o_X\leq 2i\p\bar\p|z|^2 \mbox{ in } B(z_0,2r_0) = \{|z|< 2r_0\}.
\eea
Let $x_0\in X$ be a point where $\varphi$ attains its minimum. We shall show that 
\bea
\varphi(x_0)\geq -C
\eea
where $C$ is a constant depending only on $X,\o_X,p,\gamma,
\|e^{nF}\|_{L^1(\log L)^p}$, and $\|\varphi\|_{L^1(X,\o_X)}$. It is not difficult to show by a separate argument that $\|\varphi\|_{L^1(X,\o_X)}$ is bounded by a constant depending only $n,\o_X$, for all functions with ${\sup}_X\varphi=0$ with $\lambda[h_\varphi]\in\Gamma\subset\Gamma_1=\{\lambda;\lambda_1+\cdots+\lambda_n>0\}$. The desired theorem would follow.

\medskip
We introduce now the auxiliary equation. Let $\Omega=B(x_0,2r_0)$, and for each $s$ with $0<s<2r_0^2$, set
\bea
u_s(z)=\varphi(z)-\varphi(x_0)+{1\over 2}|z|^2-s
\eea
and $\Omega_s=\{z\in\Omega; u_s(z)<0\}$.
Let $\psi_{s,\ell}$ be the solution of the following Dirichlet problem, 
\bea
(i\p\bar\p\psi_{s,\ell})^n={\tau_\ell(u_s)\over A_{s,\ell}}e^{nF(z)} \omega_X^n
\quad {\rm on}\ \Omega,
\quad \psi_{s,\ell}=0\ {\rm on}\ \p\Omega,
\eea
where 
the coefficients $A_{s,\ell}$ are defined by
\bea
A_{s,\ell}=
\int_{\Omega}
\tau_\ell(-u_s)e^{nF}\o_X^n
\ \to\ 
A_s=\int_{\Omega_s}(-u_s)e^{nF}\o_X^n,\mbox{ as }\ell \to \infty.
\eea
By \cite{CKNS}, the solution $\psi_{s,\ell}$ of this Dirichlet problem exists and is unique.

\medskip
We can now state the key comparison inequality between $u_s$ and $\psi_{s,\ell}$
\bea
\label{us:psi}
-u_s
\leq C(n,\gamma)\, A_{s,\ell}^{1\over n+1}(-\psi_{s,\ell})^{n\over n+1}
\eea
on $\bar\Omega$, where $C(n,\gamma)$ depends only on $n$ and $\gamma$. Note that the auxiliary Monge-Amp\`ere equation is of the general form considered in (\ref{eq:MA}), with $\o=0$,
$\tilde q(z)=0$, $\Lambda=0$, 
and $q(z)=\varphi(z_0)-{1\over 2}|z|^2$, (and $s\to -s$). In particular $q_{\bar kj}=-{1\over 2}\delta_{kj}$. It is then easy to apply the maximum principle to the function
\bea
\Phi=-\e(-\psi_{s,\ell})^{n\over n+1}-u_s
\eea
with the preliminary computations as in 
Lemma \ref{lm:maximum}, and the desired inequality follows with $\e^{n+1}=A_{s,\ell}\gamma^{-1}(n+1)^nn^{-2n}$.

\medskip
We can now follow the template provided by the compact K\"ahler case. First we establish the inequality
\bea
\label{exponential}
\int_{\Omega_s}{\rm exp}\{\beta {(-u_s)^{n+1\over n}\over A_s^{1/n}}\}\o_X^n
\leq C
\eea
where $\beta$ and $C$ are strictly positive constants depending only on $n,\gamma,r_0$. This is the analogue in the Hermitian case of the inequality (\ref{As1}). For this, we apply Young's inequality as in the compact K\"ahler case, with the $\alpha$-invariant replaced by the following inequality of Kolodziej \cite{K} for plurisubharmonic functions $\psi$
\bea
\int_De^{-\alpha\psi}dV\leq C
\eea 
on bounded pseudoconvex domains $D\subset{\bf C}^n$ with $\psi=0$ on $\p D$, and Monge-Amp\`ere measure $\int_D(i\p\bar\p\psi)^n=1$. Here $\alpha$ is a strictly positive constant\footnote{Kolodziej's inequality is a generalization to all dimensions of a classic inequality in one complex dimension due to Brezis and Merle \cite{BM}. Kolodziej's original proof used pluripotential theory. A more recent PDE proof has been provided by Wang, Wang, and Zhou \cite{WWZ}.}.

\medskip
Next, we can apply Young's inequality to the exponential inequality (\ref{exponential}) and obtain for any $p>n$ the following reverse H\"older inequality
\bea
A_s\leq C_0(\int_{\Omega_s}e^{nF}\o_X^n)^{1+\delta_0}
\eea
with $\delta_0={1\over n}-{1\over p}$, and $C$ is a constant depending only on $\o_X, n,p,\gamma$ and $\|e^{nF}\|_{L^1(\log L)^p}$. Setting
\bea
\phi(s)
=
\int_{\Omega_s}e^{nF}\o_X^n
\eea
we can rewrite the reverse inequality as $A_s\leq C_0(\phi(s))^{1+\delta_0}$. And since it is easily seen that $A_s
\geq t\int_{\Omega_{s-t}}e^{nF}\o_X^n$, we conclude that for any $t\in (0,s)$
\bea
\label{DGinc}
t\phi(s-t)\leq C_0(\phi(s))^{1+\delta_0}.
\eea
We can now apply a version of De Giorgi's lemma for monotone increasing and positive functions $\phi(s)$ on $(0,s_0)$ satisfying ${\rm lim}_{s\to 0^+}\phi(s)=0$ and the growth rate
(\ref{DGinc}), which says that under these conditions, there must exist a constant $c_0>0$ depending only on $s_0, C_0$ and $\delta_0$ so that 
\bea
\phi(s_0)\geq c_0.
\eea

Finally, we note the elementary inequality
\bea
\phi(s_0)\log(-\varphi(x_0)-s_0)^{1\over 2}
\leq
\int_{\Omega_{s_0}}\log{-\varphi-2^{-1}|z|^2\over (-\varphi(x_0)-s_0)^{1\over 2}}e^{nF}\o_X^n
\eea
which results itself from the elementary inequality $-\varphi-{1\over 2}|z|^2>-\varphi(x_0)-s_0$ on $\Omega_0$ since we may assume without loss of generality that $s_0<{1\over 2}$ and $\varphi(x_0)<-2$. Applying the Young's inequality to the integrand on the right hand side, we find
\bea
\phi(s_0)\log(-\varphi(x_0)-s_0)^{1\over 2}
\leq
{1\over 2}
{\rm max}\big\{
\|e^{nF}\|_{L^1(\log L)^p}+C'_pC_1(\o_X),C_p''
{\|\varphi\|_{L^1}+C_2(\o)\over
(-\varphi(x_0)-s_0)^{1\over 2}}\big\}
\nonumber
\eea
Since $\phi(s_0)$ is bounded from below by a constant, it follows that $-\varphi(x_0)-s_0$ must be bounded from above. The theorem is proved.

\section{Application to $(n-1)$-form Monge-Amp\`ere equations}
\setcounter{equation}{0}

A new generation of problems in complex geometry, notably from non-K\"ahler geometry and mathematical physics, has led to non-linear equations of the form $f(\lambda)=e^F$, but where $\lambda$ are not the eigenvalues of a hessian matrix $h_\varphi=\o_X^{-1}(\o+i\p\bar\p\varphi)$, but of a more general matrix involving the derivatives of $\varphi$ to second order. One basic example of such an equation is the so-called $(n-1)$-form Monge-Amp\`ere equation, solved by G. Sz\'ekelyhidi, V. Tosatti, and B. Weinkove \cite{STW} in their solution of the Gauduchon conjecture. Central to the study of this equation is the $L^\infty$ estimate. In this section, we show how our method of auxiliary Monge-Amp\`ere equations can be used to establish this estimate, with in fact slightly weaker and more general hypotheses than in \cite{STW}.

\medskip
Indeed, let $(X,\o_X)$ be a compact Hermitian manifold. If $\o$ is any other Hermitian metric, we consider the smooth functions $\varphi$ for which the form
\bea
\tilde\o_\varphi
=
\o+{1\over n-1}(\Delta_{\o_X}\varphi \o_X-i\p\bar\p \varphi)
\eea 
is positive, where $\Delta_{\o_X}\varphi=n{i\p\bar\p\varphi\wedge\o_X^{n-1}\over\o_X^n}$ is the rough Laplacian with respect to the metric $\o_X$. Then

\begin{theorem}
\label{tildeomega}
Let $\tilde h_\varphi=\o_X^{-1}\tilde\o_\varphi$, and $\lambda[\tilde h_\varphi]$ be its eigenvalues. Fix any $p>n$. Then any smooth solution $\varphi$ of the equation
\bea
\prod_{j=1}^n(\tilde\lambda_\varphi)_j=e^{nF}
\eea
with $\tilde\o_{\varphi}>0$ and ${\rm sup}_X\varphi=0$ satisfies
\bea
{\rm sup}_X|\varphi|\leq C
\eea
where $C$ is a constant depending only on $\o_X,n,p,\o$ and $\|e^{nF}\|_{L^1(\log L)^p(X,\o_X)}$.
\end{theorem}

We observe that this result is more precise than in \cite{TW1,TW2} in the sense that only the norm $\|e^{nF}\|_{L^1(\log L)^p(X,\o_X)}$ is needed, and not pointwise norms. In \cite{Sz}, G. Sz\'ekelyhidi obtains an estimate which depends on $\| e^{nF}\|_{L^q}$ for $q>2$, using Blocki's approach \cite{B} based on the Alexandrov-Bakelman-Pucci (ABP) maximum principle.

\smallskip

We sketch the proof. A first observation is that the tensor $\Theta^{i\bar j}$, as defined in \cite{TW1} by
\bea
\Theta^{k\bar j}
=
{1\over n-1}(({\rm Tr}_{\tilde \o_{\varphi}}\o_X) \o_X^{k\bar j}-\tilde\o_{\varphi}^{k\bar j})
\eea
is positive definite. In fact, in coordinates where $(\o_X)_{j\bar k}=\delta_{jk}$, $(\tilde \o_\varphi)_{j\bar k}=\lambda_j\delta_{jk}$ at a given point, then $\Theta^{k\bar j}={1\over n-1}(\sum_{\ell\not=j}{1\over\lambda_\ell})\delta_{jk}$.

Next, fix $r_0$ small enough so that, at any point $z_0\in X$, there is a holomorphic coordinate system $z$ with
\bea
{1\over 2}i\p\bar\p|z|^2\leq \o_X\leq 2 i\p\bar\p|z|^2,
\qquad z\in B(z_0,2r_0) = \{|z|<2r_0\}.
\eea
Let $x_0$ be a minimum point for $\varphi$, $\Omega=B(x_0,2r_0)$, and fix a small constant $\e'>0$ depending only on $X,\o_X,\o$ so that
\bea
\o\geq{2\e'\over n-1}({\rm Tr}_{\omega_X}\omega)\o_X.
\eea
Define now for $s\in (0,s_0)$, $s_0=4\e'r_0^2$,
\bea
u_s(z)=\varphi(z)-\varphi(x_0)+\e'|z|^2-s,
\qquad z\in \Omega.
\eea
Then $u_s(z)>0$ for $z\in\p\Omega$, and hence the sublevel set $\Omega_s=\{z|u_s(z)<0\}\cap \Omega$ is relatively compact in $\Omega$ and also an open set. We can now consider the solution $\psi_{s,\ell}$ of the following auxiliary Dirichlet problem on $\Omega$ for the complex Monge-Amp\`ere equation,
\bea
(i\p\bar\p\psi_{s,\ell})^n
=
{\tau_\ell(-u_s)\over A_{s,\ell}}e^{nF}\o_X^n
\quad {\rm on}\ \Omega,
\quad
\psi_{s,\ell}=0\ {\rm on}\ \p\Omega
\eea
with $A_{s,\ell}$ defined by
\bea
A_{s,\ell}=\int_{\Omega}\tau_\ell(-u_s)e^{nF}\o_X^n
\to
A_s=\int_{\Omega_s}(-u_s)e^{nF}\o_X^n.
\eea
By \cite{CKNS}, this Dirichlet problem admits a unique solution $\psi_{s,\ell}$. From the very equation, the solution $\psi_{s,\ell}$ satisfies $\int_\Omega(i\p\bar\p\psi_{s,\ell})^n=1$.

\smallskip
The key comparison estimate is now
\bea
\label{comparison2}
-u_s\leq ({n+1\over n})^{n\over n+1}A_{s,\ell}^{1\over n+1}
(-\psi_{s,\ell})^{n\over n+1}
\eea
which follows from the non-positivity on $\Omega$ of the test function
\bea
\Phi=-\e(-\psi_{s,\ell})^{n\over n+1}-u_s,
\qquad \e=A_{s,\ell}^{1\over n+1}({n+1\over n})^{n\over n+1}.
\eea 
To show that $\Phi\leq 0$, we apply the maximum principle at a maximum point $x_0$ for $\Psi$, but with the operator 
$Lv=\Theta^{j\bar k}v_{\bar kj}$, so that $L\Phi(x_0)\leq 0$. The desired inequality follows using the arithmetic-geometric inequality, along the same lines as Lemmas.

\smallskip
Once we have the key comparison estimate (\ref{comparison2}), the proof follows the same template as in the compact K\"ahler case, or more precisely the compact Hermitian case. Thus we apply the exponential estimate of Kolodziej to arrive at a reverse H\"older inequality, which allows an application of De Giorgi's lemma, and the proof is completed using a general $L^1$ bound for $-\varphi$. This completes the proof of Theorem \ref{tildeomega}.

\medskip
Finally, we note that there has been considerable interest recently in more general equations involving $f(\lambda[\tilde h_\varphi])$, motivated in part by non-K\"ahler geometry, mirror symmetry, and equations from string theories (see e.g. 
\cite{FY, CJY, Po, P}). In particular, subsolutions and general $C^2$ estimates have been obtained in \cite{G, Sz}. A frequent feature of these new equations is the appearance of gradient terms. For related developments, see \cite{FWW, GN, CHZ, PPZ0, PPZ1, PPZ} and references therein.

\section{Application to Green's functions}
\setcounter{equation}{0}

Perhaps surprisingly, the auxiliary Monge-Amp\`ere equations which we have been discussing turn out to very effective in many seemingly unrelated problems in geometric analysis. In this section, we show how they can be applied to derive lower bounds for the Green's function. The Green's function is the solution of a linear partial differential equation, so this would be a case of a linear problem solved by a non-linear method. 
Such methods had been instrumental in the study of Schr\"odinger equations \cite{CLY, SY, LY, SWYY}, and in the celebrated work of M. Kuranishi \cite{Ku1, Ku2, Ku3} on the embeddability of strongly pseudo-convex manifolds, as pointed out by C. Fefferman. However, it does not appear that comparisons with the Monge-Amp\`ere equation had been used before in linear problems. In our case, this non-linear method is needed because it gives bounds which are uniform in the underlying geometry.

\medskip
In this section, we describe some of the results obtained in \cite{GPS}. The setting is then a compact K\"ahler manifold $(X,\o_X)$ of dimension $n$, and a closed nonnegative $(1,1)$-form $\chi$ whose cohomology class $[\chi]$ is nef and big, that is, $[\chi]$ lies in the closure of the K\"ahler cone and $\int_X\chi^n>0$. Then $[\chi+t\o_X]$ is a K\"ahler class for any $t\in(0,1]$. For any $\e>0$, $N>0$, $\gamma\geq 1$, and any fixed $t\in (0,1]$, we introduce the following class of K\"ahler metrics,
\bea
{\cal M}_t'(N,\e,\gamma)
=
\{\o\in [\chi+t\o_X];{1\over V}\int_Xe^{(1+\e)F_\o}\o_X^n
\leq N,
\ 
{\rm sup}_Xe^{-F_\o}\leq \gamma^{-1}\}.
\eea
where $F_\o$ denotes the relative volume form of $\o$,
\bea\label{eqn:10.2}
F_{\o}=\log({\o^n/[\o^n]\over \o_X^n/[\o_X^n]}).
\eea
We are interested in bounds which hold uniformly for $\o\in{\cal M}'_t(N,\e,\gamma)$.
The following is a basic example of an estimate which would be well-known and readily established for a fixed K\"ahler form $\o$, but requires a non-linear method in order to have uniform bounds with respect to $\omega\in {\cal M}_t'(N,\e,\gamma)$:

\medskip
\begin{theorem} 
\label{Laplacian}
Let $\o$ be any K\"ahler form in the class
${\cal M}_t'(N,\e,\gamma)$.

{\rm (a)} Let $v\in L^1(X,\o^n)$ be any function satisfying $\int_X v\o^n=0$. Let $\Omega_s=\{v>s\}$ be the super-level sets of $v$. Assume that
\bea
v\in C^2(\bar \Omega_{0}), \qquad
\Delta_{\o_t}v\geq -a \ \ {\rm in}\ \Omega_0
\eea
for some $a>0$. Then
\bea
{\rm sup}_Xv\leq C(a+\|v\|_{L^1(X,\o^n)})
\eea
where $C$ is a constant depending only on $n,\o_X,\chi, \e,N$, and $\gamma$.

{\rm (b)} Assume now that $v\in C^2(X)$ and that
\bea
|\Delta_\o v|\leq 1, \quad{\rm and}
\quad\int_Xv\o^n=0.
\eea
Then there is a constant $C$ depending only on $n,\o_X,\chi, \e,N,\gamma$ such that $\|v\|_{L^1(X,\o^n)}\leq C$.
\end{theorem}


The proof relies again on an auxiliary complex Monge-Amp\`ere equation. Replacing $v$ by $v/a$ we may assume $a = 1$.  Thus let as before
$u\to\tau_\ell(u)$ be again a sequence of smooth strictly positive functions which decrease to the function $u\to u\chi_{{\bf R}_+}(u)$ as $\ell\to\infty$, and consider the solution $\psi_{s,\ell}$ of the following equation
\bea
(\o+i\p\bar\p\psi_{s,\ell})^n={\tau_\ell(v-s)\over A_{s,\ell}}\o^n,
\qquad
{\rm sup}_X\psi_{s,\ell}=0,
\eea
where $A_{s,\ell}$ is again a normalizing constant
\bea
A_{s,\ell}={1\over [\o^n]}\int_X\tau_\ell(v-s)\o^n
\ \to\ {1\over [\o^n]}\int_{\Omega_s}(v-s)\o^n
=A_s
\eea
as $\ell\to\infty$. Yau's theorem insures then the existence of a unique solution $\psi_{s,\ell}$.

\smallskip
We need to compare $v$ to $\psi_{s,\ell}$. For this, we express $\o$ by the $\p\bar\p$-lemma as
\bea\label{eqn:10.8}
\o=\chi+t\o_X+i\p\bar\p\varphi,\qquad {\rm sup}_X\varphi=0
\eea
for a unique function $\varphi$. Then the key estimate is the following: there exists a constant $\Lambda$ depending only on $n,p,\chi,\o_X,{\rm Ent}_p(\o)$ so that
\bea
-\psi+\varphi+\Lambda\geq 1,
\quad{\rm and}\quad\Phi\leq 0
\eea
and the function $\Phi$ is defined by
\bea
\Phi=-\e(-\psi+\varphi+\Lambda)^{n\over n+1}+v-s
\eea
with $\e^{n+1}=({n+1\over n^2})^n(a+\e n)^n A_{s,\ell}$. The proof of this is analogous to the ones we have used for the previous auxiliary complex Monge-Amp\`ere equations.

\smallskip
Once we know that $\Phi\leq 0$, we can apply an $\alpha$-invariant inequality, uniform for all K\"ahler classes bounded by a fixed multiple of $\o_X$ to obtain an inequality of the form
\bea
\int_{\Omega_s}
{\rm exp}(\alpha{(v-s)^{n+1\over n}\over A_{s,\ell}^{1/ n}})\o_X^n
\leq C
\eea
and hence, using Young's inequality, the following reverse H\"older inequality for $p>n$,
\bea
\int_{\Omega_s}(v-s)^{(n+1)p\over n}e^F\o_X^n\leq C\,A_{s,\ell}^{p\over n}\to C\,A_s^{p\over n}.
\eea
This readily implies $A_s\leq (\int_{\Omega_s}e^F\o_X^n)^{1+n\over np'}$, with $p'$ the dual exponent of $p(n+1)/n$. An easy consequence is the growth inequality
\bea
r\phi(s+r)\leq C\phi(s)^{1+\delta_0},
\qquad s\geq 0,\ r>0
\eea
for the monotone decreasing function $\phi(s)=\int_{\Omega_s}e^F\o_X^n$. Again by a De Giorgi lemma, it follows that $\phi(s)$ must vanish for $s>S_0$, where $S_0$ can be estimated by the constants in the growth condition. Thus $v\leq S_0$, and part (a) of Theorem \ref{Laplacian} is proved.

\medskip
Next, we sketch the proof of Part (b) of Theorem \ref{Laplacian}. We observe that the proof of Part (a) did not require a uniform lower bound $\gamma$ for the volume form in the definition of the class ${\cal M}'_t(N,\e,\gamma)$, but the proof of Part (b) will. 

We argue by contradiction. Thus assume that there exists a sequence of metrics $\o_j\in {\cal M}_{t_j}'(N,\e,\gamma)$ with $\{t_j\}_j\subset (0,1)$ and a sequence of functions $\hat v_j\in C^2(X)$ satisfying $\|\hat v_j\|_{L^1(X,\o_j^n)}=1$ and 
\bea
\Delta_{\o_j}\hat v_j=\hat h_j,
\quad \int_X\hat v_j\o_j^n=0,
\quad
{\rm sup}_X|\hat h_j|\to 0
\eea
as $j\to\infty$. Multiplying the above equation by $\hat v_j$ and integrating by parts gives
\bea
\int_X|\na \hat v_j|_{\o_j}^2\o_j^n
=|\int_X \hat h_j\hat v_j\o_j^n|
\leq V_{\o_j}^{1\over 2}\,{\rm sup}_X|\hat h_j|\to 0.
\eea
On the other hand, we can write
\bea
\int_X|\na\hat v_j|_{\o_X}\o_X^n
&\leq&
\int_X(|\na\hat v_j|_{\o_j}^2{\rm tr}_{\o_X}\o_j)^{1\over 2}\o_X^n
\nonumber\\
&\leq &
(\int_X|\na\hat v_j|_{\o_j}^2e^{F_j}\o_X^n)^{1\over 2}
(\int_X({\rm tr}_{\o_X}\o_j)e^{-F_j}\o_X^n)^{1\over 2}.
\eea
The first factor in the right hand side tends to $0$ as $j\to\infty$, since 
\bea
\int_X|\na\hat v_j|_{\o_j}^2e^{F_j}\o_X^n
=
{V_X\over V_{\o_j}}\,\int_X|\na \hat v_j|_{\o_j}^2\o_j^n\to 0
\eea
since $V_{\o_j}\geq [\chi^n]$ for all $\o_j$. In view of the lower bound for the volume form of metrics $\o_j\in{\cal M}_{t_j}'(N,\e,\gamma)$, we have
\bea
\int_X({\rm tr}_{\o_X}\o_j)e^{-F_j}\o_X^n
\leq
{1\over n\gamma}\int_X \o_j\wedge \o_X^{n-1}\leq
C
\eea
where $C$ is a constant independent of $j$. This implies that
\bea
\int_X|\na \hat v_j|_{\o_X}\o_X^n\to 0
\quad{\rm as}\ j\to\infty.
\eea
From this, it is not difficult to deduce that $\hat v_j$ is uniformly bounded in the Sobolev space $W^{1,1}(X,\o_X)$, and that, ater passing to subsequences, it must converge to a constant $\hat v_\infty$ in $L^1(X,\o_X^n)$. This constant can on one hand be verified to be nonzero in view of the normalization $\|\hat v_j\|_{L^1(X,\o_j^n)}=1$, and on the other hand to be $0$, in view of the condition $\int_X\hat v_j\o_j^n=0$ satisfied by $\hat v_j$. This is a contradiction and Theorem \ref{Laplacian} is proved.

\bigskip
We can now establish uniform estimates for the Green's function.
Recall now that the Green's function with respect to the K\"ahler metric $\o$ is defined by the equations
\bea
\Delta_\omega G(x,y)=-\delta_x(y)+{1\over V_\omega},
\qquad
\int_XG(x,y)\o^n=0.
\eea

\begin{theorem}
\label{Green}
Fix $\e>0, N>0$ and $\gamma\in (0,1)$. Then for any K\"ahler metric $\o\in {\cal M}_t'(\e,N,\gamma)$, the corresponding Green's function $G(x,y)$ satisfies the following estimates, with constants $C$ which depend only on $n,\o_X,\chi, \e,N$ and $\gamma$:

{\rm (a)} $\int_X|G(x,\cdot)|\o^n\leq C$;

{\rm (b)} ${\rm inf}_{y\in X} G(x,y)\geq -C$;

{\rm (c)} For any $\delta>0$, there is a constant $C_\delta$ depending additionally on $\delta$ so that, for all $x\in X$,
\bea
\int_X|G(x,\cdot)|^{{n\over n-1}-\delta}\o^n+
\int_X|\na G(x,\cdot)|^{{2n\over 2n-1}-\delta}\o^n\leq C_\delta.
\eea
\end{theorem}

\medskip
Parts (a) and (b) of Theorem \ref{Green} are direct consequences of Theorem \ref{Laplacian}. For example,
for any fixed $x$, the function $v=-G(x,y)$ is in $C^\infty(X\setminus \{x\})$ and satisfies the conditions in Part (a) of Theorem \ref{Laplacian}  with $\Delta_\o v(y)=-{1\over V_\o}\geq -{1\over [\chi^n]}$ for $y\in \{v\geq 0\}$. Thus Theorem 
\ref{Laplacian} applies, giving the lower bound 
\bea
{\rm inf}_{y\in X}G(x,y)\geq - C(1+\|G(x,\cdot)\|_{L^1(X,\o^n)})\eea
for some constant depending only on $n,p,\o_X,\chi,N$. With a bit more work, we can deduce from Part (b) of Theorem \ref{Laplacian} that $\|G(x,\cdot)\|_{L^1(X,\o^n)}\leq C$. 
Combined with the preceding inequality, we obtain Parts (a) and (b) of Theorem \ref{Green}.

\medskip
The proof of Part (c) is harder and requires a new idea, involving comparisons with another auxiliary complex Monge-Amp\`ere equation. The key inequality to be established is a uniform bound for the $L^q$ norm of $G(x,\cdot)$, 
\bea
\int_X|G(x,y)|^q\o^n(y)\leq C_q
\eea
first for $q\in 1+{1\over r_0}$ for some $r_0>n$, and then iteratively for any $q<{n\over n-1}$.
By Part (b), we can add a uniform constant to $G(x,\cdot)$ to obtain a function ${\cal G}(x,\cdot)\geq 1$. Fix $r_0>n$ and a large $k>>1$ and set
\bea
H_k(y)={\rm min}\{{\cal G}(x,y),k\}
\eea
which we assume is smooth, by smoothing it out if necessary. 
The above integral is closely related to the following integral
\bea
\int_X {\cal G}(x,y)H_k(y)^{1\over r_0}\o^n(y)
\eea
which is itself closely related to the solution $u_k$ of the following Laplace equation
\bea
\Delta_\o u_k=-H_k^{1\over r_0}+{1\over V_\o}\int_X H_k^{1\over r_0}\o^n,
\qquad
{1\over V_\o}\int_X u_k\o^n=0.
\eea

To estimate $u_k$, we introduce another auxiliary complex Monge-Amp\`ere equation,
\bea
(\o+i\p\bar\p\psi_k)^n={H_k^{n\over r_0}+1\over V_\o^{-1}\int_X(H_k^{n\over r_0}+1)\o^n}\o^n,
\qquad
{\rm sup}_X\psi_k=0.
\eea
It can then be shown that
\bea
{\rm sup}_X|\psi_k|\leq C
\eea
and that
\bea
\e'u_k+(\psi_k-\varphi)-{1\over V_\o}\int_X(\psi_k-\varphi)\o^n
\leq C
\eea
where $C$ and $\e'$ are uniform constants, and $\varphi$ is the potential for the K\"ahler metric $\o\in[\chi+t\o_X]$ introduced earlier in (\ref{eqn:10.8}). From here, the desired inequality follows.

\medskip
Finally, we establish integral bounds for $\na G(x,y)$. First we note the elementary inequality
\bea
\int_X{|\na_y{\cal G}(x,y)|_{\o(y)}^2\over {\cal G}(x,y)^{1+\beta}}\o^n(y)\leq {1\over\beta}
\eea
which holds for all $\beta>0$, and follows from applying Green's formula to $u(y)={\cal G}(x,y)^{-\beta}$. Next, setting
\bea
H_k(y)={\rm min}\{{|\na_y{\cal G}(x,y)|_{\o(y)}^2\over {\cal G}(x,y)^{1+\beta}},k\}
\eea
and arguing as in the estimate of $\|G(x,\cdot)\|_{L^{1+{1\over r_0}}(X,\o^n)}$, we can show that
\bea
\int_X{\cal G}(x,y)H_k(y)^{1\over r_0}\o^n(y)
\leq C.
\eea
The desired $L^s(X,\o^n)$ bound for $\na G(x,\cdot)$ ultimately follows. Q.E.D.

\bigskip
It is instructive to compare the preceding theorem with the classic result of Cheng and Li \cite{CL} on lower bounds for the Green's function in Riemannian geometry. This result asserts the existence of a uniform lower bound depending only the dimension, the diameter, the volume, and a lower bound on the Ricci curvature. In our K\"ahler setting, Theorem \ref{Green} is easily seen to imply the following theorem,
where only a lower bound on the scalar curvature, combined with an integral estimate for the volume form, suffice to give lower bounds for the Green's function which are uniform in $\o$:

\medskip
\begin{theorem}
\label{Kahler} Let $\o$ be any K\"ahler metric in $[\o_X]$. Then if $\|e^{F_\o}\|_{L^{1+\e}(X,\o_X)}\leq N$ for some $\e>0, N>0$, and the scalar curvature $R(\o)$ satisfies $R(\o)\geq -\kappa$ for some $\kappa\geq 0$, then the Green's function of $(X,\o)$ satisfies 
\bea
{\rm inf}_{y\in X}G(x,y)\geq -C,\quad x\in X,
\eea
for a constant $C$ depending only on $n,\o_X,\e,N$, and $\kappa$.
\end{theorem}

\medskip
We note that the class ${\cal M}_t'(\e,N,\gamma)$ of Hermitian metrics which we have used so far is not the only class to which the methods of this section apply. Other classes are described in \cite{GPS} as well, with a key difference being the replacement of a pointwise lower bound on $e^{F_\o}$ by an integral bound. In fact, for applications to diameter bounds to be described in the next section, it will be important to relax further the constant lower bound $\gamma$ to a non-negative function which may vanish along a closed set of Hausdorff dimension strictly less than $2n-1$. 

\bigskip
In the remaining part of this section, we describe some applications of the above bounds for the Green's function to a priori estimates for the complex Monge-Amp\`ere equation. A priori estimates are often obtained by applying the maximum principle to elliptic differential inequalities satisfied by the quantity under consideration. Sharp lower bounds for the Green's function can provide a more effective tool, especially if we consider conditions involving integrals. The following is a sample of sharp estimates which can be obtained in this manner.

\medskip
Recall that we have assumed that the $(1,1)$-form $\chi$ is non-negative, and the K\"ahler class $[\chi]$ is big. By Kodaira's lemma, there is an effective divisor $D$ in $X$ such that
\bea
\chi-\e_0 {\rm Ric}(h_D)\geq \delta_0 \,\o_X
\eea
for some positive constants $\e_0,\delta_0$, and a Hermitian metric $h_D$ on the line bundle $[D]$ associated with $D$. Let $s_D$ be a holomorphic section defining $D$ with ${\rm sup}_X|s_D|_{h_D}^2=1$. Let $\o$ be any metric in $[\chi+t\o_X]$. Let $\varphi$ be its potential, i.e. 
$\o=\chi+t\o_X+i\p\bar\p\varphi$ and ${\rm sup}_X\varphi=0$. Thus, in the notation (\ref{eqn:10.2}) for the relative volume form $F_\o$, $\varphi$ satisfies the complex Monge-Amp\`ere equation
\bea
(\chi+t\o_X+i\p\bar\p\varphi)^n
=
c e^{F_\o}\o_X^n,
\qquad
c={V_\o\over V_X},
\quad {\rm sup}_X\varphi=0.
\eea

\begin{theorem}
\label{C1}
Fix $\e, N, \gamma\in (0,1)$ and $p>n$. Then for any $t\in (0,1]$ and $\o\in {\cal M}_t'$, we have the estimate
\bea
|\na\varphi|_{\o_X}^2\leq {C\over |s_D|_{h_D}^{2A}}
\eea
where $C$ depends only on $n,\e,\chi,\o_X,N,\gamma,p$, and $A>0$ depends only on $n,\e,\chi,\o_X,N,\gamma$.
\end{theorem}

Note that gradient bounds had been obtained before, but under pointwise assumptions on $|\na F|$ \cite{Bl, PS09}. An earlier result requiring an $L^p$ bound for $|\na F|$ with $p\geq 2n$ is in \cite{ChHe, GPT21s}. The range $p>n$ in the above theorem is sharp.

\begin{theorem}
\label{C2}
Under the same assumptions as in Theorem \ref{C1}, but with $p>2n$, we have the estimate
\bea
|i\p\bar\p\varphi|_{\o_X}^2\leq {C\over |s_D|_{h_D}^{2B}}
\eea
where $C$ depends only on $n,\e,\chi,\o_X,N,\gamma,p$ and 
$\int_X|\na F_\o|_{\o_X}^pe^{F_\o}\o_X^n$, and $B>0$ depends only on $n,\e,\chi,\o_X,N,\gamma$.
\end{theorem}

We have formulated the above estimates for families of degenerating metrics. But even in the case of a fixed background K\"ahler form $\o$, the above estimates improve on the known ones. For example, we have

\begin{theorem}
\label{fixed}
Consider the complex Monge-Amp\`ere equation
\bea
(\o_X+i\p\bar\p\varphi)^n=e^{F}\o_X^n,
\quad {\rm sup}_X\varphi=0
\eea
on an $n$-dimensional compact K\"ahler manifold $(X,\o_X)$. Assume that $F$ satisfies the condition
\bea
&&
\|e^F\|_{L^{1+\e}(X,\o_X)}\leq N\nonumber\\
&&
{\rm sup}_X e^{-F}\leq\gamma^{-1},
\eea
Then for any $p>2n$, we have

{\rm (a)} ${\rm sup}_X|i\p\bar\p\varphi|_{\o_X}^2\leq C$, where $C>0$ is a constant depending only on $n,p,\o_X,\e,N,\gamma$ and $\int_X|\na F|_{\o_X}^pe^F\o_X^n$.

{\rm (b)} ${\rm sup}_X|\na i\p\bar\p\varphi|_{\o_X}^2\leq C$, where $C>0$ is a constant depending only on $n,p,\o_X,\e,N, \gamma$,
\bea
\gamma,\int_X|\na F|_{\o_X}^pe^F\o_X^n,
\quad
\int_X|i\p\bar\p F|_{\o_X}^pe^F\o_X^n,
\eea
and upper and lower bounds for the endomorphism $\o_X^{-1}(\o_X+i\p\bar\p\varphi)$.
\end{theorem}

For the estimate in (a) to hold in general, we do need $p\geq 2n$. We also note that previous $C^3$ bounds had required a $C^3$ bound for $F$, and that the proof of (b) also relied on the approach of \cite{PSS}, which relied on estimating the connection forms instead of the potentials.

\section{Application to diameter bounds}
\setcounter{equation}{0}

In general, estimates for the Green's function can imply estimates for the diameter of the underlying metric. This can be seen as follows \cite{GPSS}.

\medskip
Let $x_0,y_0$ be points with $d_\o(x_0,y_0)={\rm diam} (X,\o)$, and define the function $d$ on $X$ by $d(y)=d_\o(x_0,y)$. Then $d$ is a Lipschitz function with Lipschitz constant $1$. The Green's formula applied to $d(y)$ gives
\bea
d(x)
=
{1\over [\o^n]}\int_X d(y)\o(y)^n
+
\int_X\<\na_yG(x,y),\na d(y)\>_{\o(y)}\o(y)^n.
\eea
The fact that $d(x_0)=0$ gives 
\bea
{1\over [\o^n]}\int_X d(y)\o(y)^n
=-
\int_X\<\na_yG(x_0,y),\na d(y)\>_{\o(y)}\o(y)^n
\leq
\int_X|\na_yG(x_0,y)|_{\o(y)}\o(y)^n.
\nonumber
\eea
We can then write
\bea
{\rm diam}(X,\o)
=d(y_0)
&=&{1\over [\o^n]}\int_X d(y)\o(y)^n
+
\int_X\<\na_yG(y_0,y),\na d(y)\>_{\o(y)}\o(y)^n
\nonumber\\
&\leq &
\int_X|\na_yG(x_0,y)|_{\o(y)}\o(y)^n
+
\int_X|\na_yG(y_0,y)|_{\o(y)}\o(y)^n
\nonumber
\eea
which shows that the diameter can be estimated by an integral bound for the gradient of the Green's function.

\medskip
Thus the bounds obtained in the previous section already imply some diameter bounds. However, for many geometric applications, such as diameters in the K\"ahler-Ricci flow, it is important to relax the conditions on the lower bound $\gamma$ for the volume form. It turns out that this is possible, albeit quite non-trivial. Thus the following theorems were established in \cite{GPSS}:

\begin{theorem}
\label{diameter}
Let $(X,\o_X)$ be an $n$-dimensional connected compact K\"ahler manifold. For given parameters $A, K>0, p>n$ and $\gamma$ a continuous non-negative function, let ${\cal V}(X,\o_X,n,$ $A,p,K,\gamma)$ be the following space of K\"ahler metrics
\bea
{\cal V}(X,\o_X,n,A,p,K,\gamma)=\{\o; [\o]\cdot[\o_X]^{n-1}\leq A,\,
{\cal N}_{X,\o_X,p}(\o)\leq K,\,
{\o^n\over\o_X^n}\geq \gamma V_\o\}
\eea
where ${\cal N}_{X,\o_X,p}(\o)$ is the $p$-Nash entropy, defined by
\bea
\label{Nash}
{\cal N}_{X,\o_X,p}(\o)
=
\int_X|F|^pe^F\o_X^n
=
\|e^F\|_{L^1(\log L)^p(\o_X)},
\quad 
F= {1\over V_\o}{\o^n\over\o_X^n}.
\eea
Assume that
\bea
{\rm dim}_{\cal H}\{\gamma=0\}<2n-1
\eea
where ${\rm dim}_{\cal H}$ denotes the Hausdorff dimension. Then for any $A,K>0$ and $p>n$, there exist constants $C,c>0$ depending only on $X,\o_X,n,A,p,K,\gamma$ and $\alpha$ depending only on $n$ and $p$ so that

{\rm (a)} $\int_X|G(x,\cdot)|\o^n
+\int_X|\na G(x,\cdot)|\o^n+(-{\rm inf}_{y\in X}G(x,y))V_\o
\leq C$;

{\rm (b)} ${\rm diam}\,(X,\o)\leq C$;

{\rm (c)} ${Vol_\o (B_\o(x,R))\over Vol_\o(X)}\geq c\,R^\alpha$ for any $x\in X$, $R\in (0,1)$.
\end{theorem}

We stress that the above theorem can give bounds on the diameter even when no lower bound on the Ricci curvature is available. This is of particular importance for the K\"ahler-Ricci flow. More generally, we obtain the following K\"ahler analogue of Gromov's precompactness theorem for metric spaces:

\begin{theorem}
Let $(X,\o_X)$ be a connected $n$-dimensional K\"ahler manifold, and let $\gamma$ be a non-negative function with
\bea
{\rm dim}_{\cal H}\{\gamma=0\}<2n-1.
\eea
Then for any $A,K>0$, $p>n$, any sequence $\{\o_j\}$ in 
${\cal V}(X,\o_X,n,A,p,K,\gamma)$ admits a subsequence converging in Gromov-Hausdorff topology to a compact metric space $(X_\infty, d_\infty)$.
\end{theorem}

Several applications of these theorems to the K\"ahler-Ricci flow and to the asymptotic behavior of fibrations near the singular fibre can be found in \cite{GPSS}.

\newcommand{\neqref}[1]{(\ref{#1})}
\newcommand{\tg}{\tilde g}
\newcommand{\beqref}[1]{(\ref{#1})}

\section{Application to the taming of symplectic forms}
\setcounter{equation}{0}

Let $(X,J)$ be a compact {\em almost complex} manifold with $J$ the {\em almost complex structure}. Suppose $m = 2n$ is the real dimension of $X$. A Riemannian metric $\tilde g$ on $X$ is called {\em almost K\"ahler}, if $\tilde g$ is $J$-compatible, i.e. $\tilde g(JY, JZ) = \tilde g(Y, Z)$ for any vector fields $Y, Z$ and the associated $2$-form $\omega_{\tilde g}$ defined by $\omega_{\tilde g}(Y,Z) =\tilde g(JY, Z)$ is {\em closed}.

Let $\Omega$ be a {\em taming symplectic form}, that is,
$\Omega(Y,JY)>0$ for $Y\not=0$,
and $g$ be the associated {\em almost Hermitian metric} of $\Omega$, i.e.
$$g(Y,Z) = \frac{1}{2}\Omega(Y, JZ) + \frac{1}{2} \Omega(Z, JY), \quad \forall \mbox{ vector fields }Y, Z.$$
Write $dV_g$ for the volume form of the Riemannian metric $g$. For a smooth function $F\in C^\infty(X)$ normalized by $\int_X e^{F} dV_g = \int_X \Omega^n$, we consider the following Calabi-Yau equation on $X$ 
\begin{equation}\label{eqn:MA}
\det \tilde g = e^{2F} \det g,
\end{equation}
where we require that $\omega_{\tilde g}$ is an {\em almost K\"ahler form} with $[\omega_{\tilde g}] = [\Omega]$. As shown by Donaldson \cite{D}, the existence of solutions to this equation would have important consequences in symplectic geometry.

\medskip

It has been proved by Tosatti, Weinkove, and Yau in \cite{TWY} that the $C^2$-{\em{}a priori} estimates of $\tilde g$ satisfying the equation \neqref{eqn:MA} can be derived by the $L^\infty$ estimates of $\varphi\in C^\infty(X)$, which satisfies the linear equation 
\begin{equation}\label{eqn:varphi}
\Delta_{\tilde g} \varphi = 2n - 2n \frac{\omega_{\tilde g}^{n-1}\wedge \Omega}{\omega_{\tilde g}^n} = 2n - \tr_{\tilde g} g,\quad {\rm sup}_X \varphi = 0.
\end{equation} 
Here $\Delta_{\tilde g}$ is the usual Riemannian Laplacian operator of the metric $\tilde g$. Note that the equation \neqref{eqn:varphi} admits a unique solution since the function on the right-hand-side has integral zero against the volume form $dV_{\tilde g}$.

\smallskip

We assume $$\| e^{F}\|_{L^2(X, dV_g)}^2 = \int_X e^{2F} dV_g \le K$$for some constant $K>0$. Our main result is the following $L^\infty$ estimate of $\varphi$ in \neqref{eqn:varphi}:

\begin{theorem}\label{thm:mainnew}
Suppose $\tilde g$ is an {\em almost K\"ahler} metric solving the equation \neqref{eqn:MA} and $\varphi$ solves \neqref{eqn:varphi}. Then there exists a constant $C>0$ depending on $n, g$ and $K$ such that
$${\rm sup}_X |\varphi|\le C(1+\| \varphi\|_{L^1(X, e^{2F}dV_g)}).$$
\end{theorem}


In the following, we denote $\{x^1,\ldots, x^m\}$ a local {\em real} coordinates on some open subset of $X$. Let $\tilde g$ be an {\em almost K\"ahler} metric, and $\tilde \omega(Y,Z) = \tilde g(JY, Z)$ be the associated symplectic $2$-form. Locally we have
\begin{equation}\label{eqn:2.1}
\tilde g = \tilde g_{ij} dx^i\otimes dx^j,\quad \tilde \omega = \frac 1 2\tilde\omega_{ij}dx^i\wedge dx^j,\quad J = J_i^j dx^i\otimes \frac{\partial }{\partial x^j},
\end{equation}
where the summations are taken over $i, j\in \{1,2,\cdots, m\}$. It follows from straightforward calculations that
\begin{equation}\label{eqn:2.2}
\tilde \omega_{ij} = \tilde g_{ik} J_j^k,\quad \tilde \omega_{ij} = -\tilde \omega_{ji}.
\end{equation}
$\tilde \omega$ being {\em almost K\"ahler} means $$0= d\tilde \omega  = \frac{1}{2} \frac{\partial \tilde \omega_{ij}}{\partial x^l} dx^l\wedge dx^i \wedge dx^j = \frac{1}{6}(\frac{\partial \tilde \omega_{ij}}{\partial x^l} + \frac{\partial \tilde \omega_{li}}{\partial x^j} + \frac{\partial \tilde \omega_{jl}}{\partial x^i} )dx^l\wedge dx^i \wedge dx^j$$
in other words, 
\begin{equation}\label{eqn:2.3}
\frac{\partial \tilde \omega_{ij}}{\partial x^l} + \frac{\partial \tilde \omega_{li}}{\partial x^j} + \frac{\partial \tilde \omega_{jl}}{\partial x^i} = 0,\quad \forall i, j, l.
\end{equation}
Multiplying both sides of \neweqref{eqn:2.3} by $\tilde \omega^{ij}: = \tilde g^{ik} J_k^j$ which is skew-symmetric in $i,j$, and taking summation over $i, j$, we get
\begin{equation}\label{eqn:2.4}
\tilde g^{ik} J_k^j \frac{\partial \tilde \omega_{ij}}{\partial x^l} + 2 \tilde g^{ik} J_k^j \frac{\partial \tilde \omega_{jl}}{\partial x^i} = 0.
\end{equation}
Substituting \neweqref{eqn:2.2} to \neweqref{eqn:2.4}, we obtain
\bea\nonumber
0&  = & \tg^{ik} J_{k}^j \frac{\partial \tg_{ip}}{\partial x^l} J_j^p  + \tg^{ik} J_k^j \tg_{ip} \frac{\partial J_j^p}{\partial x^l} + 2 \tg^{ik} J_k^j \frac{\partial \tg_{jp}}{\partial x^i} J_l^p + 2 \tg ^{ik} J_k^j \tg_{jp} \frac{\partial J_l^p}{\partial x^i}\\
& = & \nonumber
- \tg^{ik}  \frac{\partial \tg_{ik}}{\partial x^l} + J_k^j  \frac{\partial J_j^k}{\partial x^l}  + 2 \tg ^{ik} J_k^j \tg_{jp} \frac{\partial J_l^p}{\partial x^i} + 2 \tg^{ik}\frac{\partial \tg_{kl}}{\partial x^i} \\
&& \nonumber - 2 \tg^{ik}\tg_{jp} J_l^p \frac{\partial J_k^j}{\partial x^i} - 2 \tg^{ik} \tg_{jp} J_k^j \frac{\partial J_l^p}{\partial x^i}\\
& = & \nonumber
- \tg^{ik}  \frac{\partial \tg_{ik}}{\partial x^l} + 2 \tg^{ik}\frac{\partial \tg_{kl}}{\partial x^i} + J_k^j  \frac{\partial J_j^k}{\partial x^l}  -2 J_p^i \frac{\partial J_l^p}{\partial x^i}  + 2 \tg^{ik}\tg_{lp} J_j^p \frac{\partial J_k^j}{\partial x^i} + 2 J_p^i \frac{\partial J_l^p}{\partial x^i},\\& = & \nonumber
- \tg^{ik}  \frac{\partial \tg_{ik}}{\partial x^l} + 2 \tg^{ik}\frac{\partial \tg_{kl}}{\partial x^i} + J_k^j  \frac{\partial J_j^k}{\partial x^l}  + 2 \tg^{ik}\tg_{lp} J_j^p \frac{\partial J_k^j}{\partial x^i},
\eea
from which we derive that
\begin{equation}\label{eqn:2.5}
 \tg^{ik}\frac{\partial \tg_{kl}}{\partial x^i}  - \frac{1}{2} \tg^{ik}  \frac{\partial \tg_{ik}}{\partial x^l} = - \frac{1}{2}J_k^j  \frac{\partial J_j^k}{\partial x^l}  -  \tg^{ik}\tg_{lp} J_j^p \frac{\partial J_k^j}{\partial x^i}.
\end{equation}
It then follows that
\begin{equation}\label{eqn:gamma}
\tg^{ik} \tilde\Gamma_{ik}^q =  \tg^{ql} (\tg^{ik} \frac{\partial \tg_{kl}}{\partial x^i} - \frac{1}{2} \tg^{ik} \frac{\partial \tg_{ik}}{\partial x^l}    ) = - \frac{1}{2} \tg^{ql}J_k^j  \frac{\partial J_j^k}{\partial x^l}  -  \tg^{ik} J_j^q \frac{\partial J_k^j}{\partial x^i}.
\end{equation}
From \neweqref{eqn:gamma}, we see that the second term in the Laplacian of a function $u$, $\Delta_{\tilde g} u =\tilde g^{ij}\frac{\partial ^2 u}{\partial x^i \partial x^j} - \tilde g^{ik} \tilde \Gamma_{ik}^q \frac{\partial u}{\partial x^q}$ is independent of the first order derivatives of the metric coefficients $\tilde g_{ij}$.

\medskip

Let $x_0\in X$ be a minimum point of $\varphi$, i.e. $\varphi(x_0) = \inf_X \varphi$. Choose a fixed number $r_0>0$ such that $2r_0\le $ the injectivity radius of the fixed Riemannian manifold $(X,g)$. Take the normal coordinates of $(X,g)$ centered at $x_0,$ $(U, \{x^1,\ldots, x^m\})$. Without loss of generality we may assume that on $U$ the following holds
\begin{equation}\label{eqn:normal}
\frac 1 2 \delta_{ij} \le g_{ij} \le 2 \delta_{ij},
\end{equation}
and the {\em Euclidean ball} $$B(x_0, 2r_0) = \{x\in U: |x|< 2r_0\}\subset \subset U$$
where $|x| = \sqrt{\sum_{j=1}^m (x^j)^2}$ is the usual Euclidean norm of the coordinates $x = (x^1,\ldots, x^m)$. We also have a constant $C_J '>0$ depending on $J, g$ such that 
\begin{equation}\label{eqn:2.8}
{\rm sup}_U\Big( | \sum_{j,k}J_k^j  \frac{\partial J_j^k}{\partial x^l}  |_{g} + |\sum_jJ_j^q \frac{\partial J_k^j}{\partial x^i}|_g\Big)\le C_J'.
\end{equation}
From the equation \neqref{eqn:gamma}, we have for {\em any} smooth function $\psi$ on $U$, there exists a uniform constant $C_J>0$ such that on $U$
\begin{equation}\label{eqn:2.9}
 | \tg^{ik} \tilde \Gamma^q _{ik} \psi_q  |\le C_J |\nabla \psi|_g \cdot \tr_{\tg} g,
\end{equation}
where $\psi_q = \frac{\partial \psi}{\partial x^q}$ and $|\nabla \psi|_g^2 = g^{ij} \psi_i \psi_j$ is the gradient of $\psi$ with respect to the {\em fixed} metric $g$. We emphasize that the constant $C_J$ in \neweqref{eqn:2.9} depends {\em only} on $g, J$ and can be made to be independent of the choice of coordinates, though the LHS of \neqref{eqn:2.9} is only locally defined. Indeed, we can see from \neqref{eqn:2.8} that $C_J'$ depends on $|J|_g$ and $|\nabla _g J|_g$, both of which are globally defined.

\medskip

Let $\eta \in (0,1)$ be a small positive constant to be determined. For any $0< s \le s_0 = \eta r_0^2$, we consider the function defined on $B(x_0, 2r_0)$
\begin{equation}\label{eqn:us}
u_s(x):= \varphi(x) - \varphi(x_0) + \eta |x|^2 - s. 
\end{equation}
By the choice of the point $x_0$, it is clear that $u_s\ge -s$. Define the sublevel set of $u_s$ by
\begin{equation}
\label{eqn:2.11}
\Omega_s: = \{x\in B(x_0, 2r_0)| ~ u_s(x)<0\}.
\end{equation}
Note that $x_0\in \Omega_s$ so $\Omega_s$ is a {\em nonempty} open subset of $B(x_0, 2r_0)$. Moreover, by the choice of $s\le \eta r_0^2$, we see that on $B(x_0, 2r_0)\backslash B(x_0, r_0)$
$$u_s\ge \eta r_0^2 - s \ge 0,$$ 
hence we have $\Omega_s\subset B(x_0, r_0)$. 
We solve the following {\em real} Monge-Amp\`ere equation on the Euclidean ball $B(x_0, 2r_0)$
\begin{equation}\label{eqn:rMA}
\det \Big( \frac{\partial^2 \psi_{s,\ell}}{\partial x^i \partial x^j}\Big) = \frac{\tau_\ell(-u_s)}{A_{s,\ell}} e^{2F} \det g,\quad \mbox{in }B(x_0, 2r_0),
\end{equation}
$\psi_{s,\ell} = 0$ on $\partial B(x_0, 2r_0)$. Here $\psi_{s,\ell}$ is strictly convex in the ball $B(x_0, 2r_0)$ and 
$$A_{s,\ell} = \int_{B(x_0, 2r_0)} \tau_\ell(-u_s) e^{2F} (\det g ) dx>0$$
is chosen so that $ \int_{B(x_0, 2r_0)}  \det\Big( \frac{\partial^2 \psi_{s,\ell}}{\partial x^i \partial x^j}\Big) dx = 1$. Note that as $\ell \to\infty$ $$A_{s,\ell}\to A_s:= \int_{\Omega_s} (-u_s) e^{2F} \det g$$
and $A_s\le C(n,g) K\le \frac 1 2C_1 $, where  $C_1>0$ depends only on $n, g$ and $K$. Hence for all $\ell$ sufficiently large which we always assume, we have \begin{equation}\label{eqn:A large}A_{s,\ell}\le C_1.\end{equation}
Let $\beta_n$ be the volume of the unit ball in ${\mathbb R}^{2n}$.
\begin{lemma}\label{lemma 1}
There exist a constant $C_2=C_2(n)>0$ such that 
$$-\inf_{B(x_0, 2r_0)} \psi_{s,\ell} \le C_2 r_0,\quad {\rm sup}_{B(x_0, r_0)} |\nabla \psi_{s,\ell}|\le C_2.$$
\end{lemma}
\noindent {\em Proof.} Since by the definition of $\psi_{s,\ell}$, $ \int_{B(x_0, 2r_0)} \det \Big( \frac{\partial^2 \psi_{s,\ell}}{\partial x^i \partial x^j}\Big) dx = 1$, so it follows from the standard ABP maximum principle
$$-\inf_{B(x_0, 2r_0)} \psi_{s,\ell} \le -\inf_{\partial B(x_0, 2r_0)} \psi_{s,\ell} + \frac{4 r_0}{\beta_n} \Big[ \int_{B(x_0, 2r_0)}\det  \Big( \frac{\partial^2 \psi_{s,\ell}}{\partial x^i \partial x^j}\Big) dx \Big]^{1/m} = \frac{4}{\beta_n} r_0.$$

\smallskip
To see the second inequality, for any fixed point $x\in B(x_0, r_0)$, denote $V = \frac{D\psi_{s,\ell}(x)}{|D\psi_{s,\ell}(x)|}$ to be the unit vector in the direction of $D\psi_{s,\ell}(x)$ (if $D\psi_{s,\ell}(x) = 0$, there is nothing to prove, so here we assume $D\psi_{s,\ell}(x)\neq 0$). Then clearly $|D\psi_{s,\ell}(x)| = D\psi_{s,\ell}(x) \cdot V$. Consider the half line $L:$ $0\le t\mapsto x + t V$ which intersects $\partial B(x_0, r_0)$ and $\partial B(x_0, 2r_0)$ at $L(t_1)$, $L(t_2)$, respectively. We have $t_2 - t_1\ge r_0$ and $0>\psi_{s,\ell}(L(t_1))\ge -\frac{4}{\beta_n} r_0$ and $\psi_{s,\ell}(L(t_2)) = 0$. Then by the convexity of the function $t\mapsto \psi_{s,\ell}(L(t))$, we have 
$$|D\psi_{s,\ell}(x)| = \psi_{s,\ell}(L(t))'|_{t= 0}\le \frac{\psi_{s,\ell}(L(t_2)) - \psi_{s,\ell}(L(t_1))}{t_2 - t_1} \le \frac{4}{\beta_n}.$$
Taking supremum over all $x\in B(x_0, r_0) $ finishes the proof of the lemma. Q.E.D.

\medskip

Take positive constants
\begin{equation}\label{eqn:lambda}
\Lambda = \frac{2n}{1+2n} (10 C_J C_2)^{2n+1} A_{s,\ell}   ,\quad \varepsilon = \Big(\frac{2n+1}{2n} \Big)^{\frac{2n}{2n+1}} A_{s,\ell}^{\frac{1}{2n+1}}
\end{equation}
where $C_2>0$ is the constant in Lemma \ref{lemma 1} and $C_J>0$ as in \neweqref{eqn:2.9}. We observe that by \neweqref{eqn:A large}, it holds that $\Lambda$ is bounded above by the uniform constant $\frac{2n}{1+2n} (10 C_J C_2)^{2n+1} C_1$.

\smallskip

Define a function $\Phi$ on $B(x_0, 2r_0)$ by 
\begin{equation}\label{eqn:Phi}
\Phi(x) = - \varepsilon (-\psi_{s,\ell}(x) + \Lambda)^{\frac{2n}{2n+1}} - u_s(x),\quad \forall x\in B(x_0, 2r_0).
\end{equation}
We claim that $\Phi\le 0$ on this ball. As a continuous function, $\Phi$ achieves its maximum at some point $x_{\max}\in {\overline{B(x_0, 2r_0)}}$. If $x_{\max}\not\in \Omega_s$, then by the definition of $\Omega_s$, clearly we have $\Phi(x_{\max})<0$. So we assume $x_{\max}\in \Omega_s\subset B(x_0, r_0)$. By the maximum principle, it follows that $\frac{\partial \Phi}{\partial x^q}\Big|_{x_{\max}} = 0$ and $\frac{\partial^2\Phi}{\partial x^i \partial x^j}\Big|_{x_{\max}}\le 0$. Hence at $x_{\max}$
$$\Delta_{\tilde g} \Phi = \tg^{ij} \frac{\partial^2 \Phi}{\partial x^i \partial x^j} - \tg^{ij}\tilde \Gamma _{ij} ^q \frac{\partial \Phi}{\partial x^q} =  \tg^{ij} \frac{\partial^2 \Phi}{\partial x^i \partial x^j} \le 0.$$
We then calculate at $x_{\max}$. 
\bea
\nonumber
0 & \ge & \Delta_{\tilde g} \Phi = \frac{2n\varepsilon}{2n + 1} (-\psi_{s,\ell} + \Lambda)^{-\frac{1}{2n+1}} \Delta_{\tilde g} \psi_{s,\ell} - \Delta_{\tilde g} \varphi - \eta \Delta_{\tilde g} |x|^2 \\
&&\nonumber+ \frac{2n\varepsilon}{(2n+1)^2} (-\psi_{s,\ell} + \Lambda)^{-\frac{2n+2}{2n+1}}|\nabla \psi_{s,\ell}|^2_{\tilde g}\\
&\ge & \frac{2n\varepsilon}{2n + 1} (-\psi_{s,\ell} + \Lambda)^{-\frac{1}{2n+1}} \Delta_{\tilde g} \psi_{s,\ell} - 2n + \tr_{\tg} g - \eta \Delta_{\tilde g} |x|^2.\label{eqn:2.16}
\eea
We first look at the term $-\eta\Delta_{\tilde g}|x|^2$ in \neqref{eqn:2.16}. It satisfies
\begin{equation}
\label{eqn:2.17}
-\eta\Delta_{\tilde g}|x|^2= - 2 \eta \tilde g^{ij} \delta_{ij} +2 \eta \tg^{ij} \tilde \Gamma_{ij}^q x^q\ge -4 \eta \tr_{\tg} g - 2 \eta C_J r_0 \tr_{\tg }g\ge -\frac{1}{10} \tr_{\tilde g} g,
\end{equation}
where $C_J>0$ is the uniform constant in \neweqref{eqn:2.9} and we have chosen $\eta>0$ such that $$\eta(4+ 2 C_J r_0)= 1/10$$ and this fixes the uniform constant $\eta$. We will denote $D^2_{ij}\psi_{s,\ell} = \frac{\partial^2 \psi_{s,\ell}}{\partial x^i \partial x^j}$ to be (Euclidean) Hessian of the function $\psi_{s,\ell}$.
For the first term in \neweqref{eqn:2.16}, we have
\bea
\nonumber
&&\frac{2n\varepsilon}{2n + 1} (-\psi_{s,\ell} + \Lambda)^{-\frac{1}{2n+1}} \Delta_{\tilde g} \psi_{s,\ell}\\
&= \nonumber&
\frac{2n\varepsilon}{2n + 1} (-\psi_{s,\ell} + \Lambda)^{-\frac{1}{2n+1}} \tg^{ij} \Big( \frac{\partial^2 \psi_{s,\ell}}{\partial x^i \partial x^j} \Big) + \frac{2n\varepsilon}{2n + 1} (-\psi_{s,\ell} + \Lambda)^{-\frac{1}{2n+1}} \tg^{ij}\tilde \Gamma_{ij}^q \frac{\partial \psi_{s,\ell}}{\partial x^q}\\
&\ge \nonumber & \frac{4n^2\varepsilon}{2n + 1} (-\psi_{s,\ell} + \Lambda)^{-\frac{1}{2n+1}} \Big(\frac{ \det D^2 \psi_{s,\ell} }{\det \tilde g} \Big)^{1/2n} -  \frac{2n\varepsilon \Lambda^{-\frac{1}{2n+1}}}{2n+1} C_J | \frac{\partial \psi_{s,\ell}}{\partial x^q}  | \cdot\tr_{\tg} g \\
&\ge \nonumber & \frac{4n^2\varepsilon}{2n + 1} (-\psi_{s,\ell} + \Lambda)^{-\frac{1}{2n+1}} \Big(\frac{-u_s}{A_{s,\ell}} \Big)^{1/2n} -  \frac{2n\varepsilon \Lambda^{-\frac{1}{2n+1}}}{2n+1} C_J C_2\cdot \tr_{\tg} g\\
&\ge &\nonumber \frac{4n^2\varepsilon^{1+ \frac{1}{2n}}}{(2n + 1) A_{s,\ell}^{1/2n}}  \Big(\frac{-u_s}{ \varepsilon (-\psi_{s,\ell} + \Lambda)^{2n/(2n+1)}} \Big)^{1/2n} -  \frac{1}{10} \tr_{\tg} g\\
&= &2n  \Big(\frac{-u_s}{ \varepsilon (-\psi_{s,\ell} + \Lambda)^{2n/(2n+1)}} \Big)^{1/2n} -  \frac{1}{10} \tr_{\tg} g\label{eqn:2.18}.
\eea
Combining the equations \neweqref{eqn:2.18}, \neweqref{eqn:2.17}, \neweqref{eqn:2.16}, we see that at $x_{\max}$ 
$$0\ge 2n  \Big(\frac{-u_s}{ \varepsilon (-\psi_{s,\ell} + \Lambda)^{2n/(2n+1)}} \Big)^{1/2n} - 2n + \frac{4}{5} \tr_{\tg} g$$
from which we easily derive that $\frac{(-u_s)}{ \varepsilon (-\psi_{s,\ell} + \Lambda)^{2n/(2n+1)}} <1$, that is, $\Phi|_{x_{\max}}<0$.  Hence we finish the proof of the claim that $\Phi\le 0$. In particular on $\Omega_s$ it holds that 
\begin{equation}\label{eqn:2.19}
- u_s \le C(n) A_{s,\ell}^{\frac{1}{2n+1}} (-\psi_{s,\ell} + \Lambda)^{\frac{2n}{2n+1}}\le C_3 A_{s,\ell}^{\frac{1}{2n+1}},
\end{equation}
for some $C_3>0$ that depends on $n, g, J, K$. Here we have applied Lemma \ref{lemma 1} to see that $|\psi_{s,\ell}|\le C_2 r_0$ and the fact that $\Lambda$ is a uniformly bounded constant. Letting $\ell\to \infty$ we conclude from \neweqref{eqn:2.19} that
\begin{equation}\label{eqn:2.20}
- u_s \le C_3 A_{s}^{\frac{1}{2n+1}},
\end{equation}
Integrating both sides of \neqref{eqn:2.20} agaist the measure $e^{2F}(\det g ) dx$ over $\Omega_s$, we get 
\begin{equation}\label{eqn:2.21}
A_s = \int_{\Omega_s} (-u_s) e^{2F} (\det g) dx\le C_3 A_{s}^{\frac{1}{2n+1}} \int_{\Omega_s} e^{2F} (\det g) dx.
\end{equation}
So we have 
\begin{equation}\label{eqn:additional}A_{s} \le C_3^{\frac{2n+1}{2n}} \Big( \int_{\Omega_s} e^{2F} (\det g) dx \Big)^{1+ \frac{1}{2n}} = C_4 \phi(s)^{1+ \frac{1}{2n}},\end{equation} where we denote $\phi(s) =  \Big( \int_{\Omega_s} e^{2F} (\det g) dx \Big)^{1+ \frac{1}{2n}} $. On the other hand, for any $0< t < s$, on the open set $\Omega_{s-t}$ we have
$$u_s(x) = u_{s-t}(x)  - t< -t, \mbox{ i.e. } -u_s(x)> t$$
It is elementary that $$A_s \ge \int_{\Omega_{s-t}} (-u_s) e^{2F}(\det g) dx \ge t \phi(s-t).$$
Combining the above, we see that 
\begin{equation}\label{eqn:2.22}
t \phi(s-t) \le C_4 \phi(s)^{1+ \frac{1}{2n}},\quad \forall \,0<t < s\le s_0.
\end{equation}
It is not hard to see that $\phi(s)$ is an {\em increasing} and {\em continuous} function in $s\in (0,s_0]$ and $\phi(s)>0$ for any $s\in (0,s_0]$ and $\lim_{s\to 0^+}\phi(s) = 0$. We can apply a version of De Giorgi's lemma to show that $\frac{2 C_4}{1-2^{-1/2n}}\phi(s_0)^{1/2n} \ge s_0$. Hence there is a uniform constant $c_0>0$ such that
\begin{equation}\label{eqn:2.23}
\phi(s_0)\ge c_0>0.
\end{equation}
Applying \neqref{eqn:additional} with $s=s_0$, we obtain $A_{s_0} \le C_5$ for a constant $0<C_5 = C_4 (2^{2n}\beta_n K)^{1+ \frac{1}{2n}}$. From the definition of $A_{s_0}$, we derive that
\begin{equation}\label{eqn:2.25}
(-\varphi(x_0)) \cdot \phi(s_0) \le s_0 \phi(s_0) + \int_{\Omega_{s_0}} (-\varphi) e^{2F}( \det g) dx + C_5.
\end{equation}
The equation \neweqref{eqn:2.3} then implies that 
\begin{equation}
-{\rm inf}_X \varphi = -\varphi(x_0)\le C_6 + C_7 \int_X (-\varphi) e^{2F}dV_g,
\end{equation}
for some constant $C_7>0$. So we have proved the inequality for the solution $\varphi$ to the equation \neweqref{eqn:varphi}
\begin{equation}\label{eqn:2.26}
{\rm sup}_X |\varphi|\le C_8 (1 + \| \varphi\|_{L^1(X,e^{2F}dV_g)}).
\end{equation}
This finishes the proof of Theorem \ref{thm:mainnew}. Q.E.D.


\bigskip

\noindent Department of Mathematics \& Computer Science, Rutgers University, Newark, NJ 07102 USA

\noindent {bguo@rutgers.edu}

\medskip

\noindent Department of Mathematics, Columbia University, New York, NY 10027 USA

\noindent {phong@math.columbia.edu}


\begin{thebibliography}{99}

{\footnotesize

\bibitem{AO} S. Abja and G. Olive, ``Local regularity for concave homogeneous complex degenerate elliptic equations comparable to the Monge-Amp\`ere equation'', arXiv:2102.07553v1.

\bibitem{B} Z. Blocki, ``On the uniform estimate in the Calabi-Yau theorem II", Science China Math. 54 (2011) 1375-1377.

\bibitem{Bl}Z. Blocki, ``A gradient estimate in the Calabi-Yau theorem'', Math. Ann. 344 (2009),
no. 2, 317 - 327.


\bibitem{BM} H.  Brezis and F. Merle, ``Uniform estimates and blowup behavior for solutions of $\Delta u = V(x) e^u$ in two dimensions'', Comm. Part. Diff. Eqn. 16(1991), 1223-1253.


\bibitem{BEGZ}S. Boucksom, P. Eyssidieux, V. Guedj and A. Zeriahi, ``Monge-Amp\`ere equations in big cohomology classes'', Acta Math. 205 (2010), no. 2, 199 - 262.

\bibitem{CKNS} L. Caffarelli, J. J. Kohn, L. Nirenberg, and J. Spruck, ``The Dirichlet problem for nonlinear second-order elliptic equations. II. Complex Monge-Amp\`ere, and uniformly elliptic, equations''. Comm. Pure Appl. Math. 38 (1985), no. 2, 209 - 252.

\bibitem{CNS} L. Caffarelli, L. Nirenberg, and J. Spruck,
``The Dirichlet problem for nonlinear second-order elliptic equations. I. Monge-Amp\`ere equation", Comm. Pure Appl. Math. 37 (1984), no. 3, 369 - 402.


\bibitem{CC} X.X. Chen and J.R. Cheng, ``On the constant scalar curvature K\"ahler metrics I - a priori estimates",
J. Amer. Math. Soc. (2021) DOI: https://doi.org/10.1090/jams/967, arXiv: 1712.06697.

\bibitem{CC1} X.X. Chen and J.R. Cheng, ``The $L^\infty$ estimates for parabolic complex Monge-Amp\`ere and Hessian equations'', arXiv:2201.13339.

\bibitem{ChHe}X.X. Chen and W. He, ``The complex Monge-Amp\`ere equation on compact K\"ahler manifolds'', Math. Ann. 354 (2012), 1583 - 1600

\bibitem{CL}S. Y. Cheng and P. Li, ``Heat kernel estimates and lower bound of eigenvalues'',
Comment. Math. Helv., 56 (1981), 327-338

\bibitem{CLY} S.Y. Cheng, P. Li, and S.T. Yau,
``On the upper estimate of the heat kernel of a complete Riemannian manifold", Amer. J. Math. 103 (1981), no. 5, 1021 - 1063

\bibitem{Ch} P. Cherrier, ``\'Equations de Monge-Amp\`ere sur les vari\'et´\'s Hermitiennes compactes'', Bull. Sc. Math (2) 111 (1987), 343 - 385.

\bibitem{CHZ} J. Chu, L. Huang, and X.H. Zhu, ``The Fu-Yau equation in higher dimensions", Peking Math. J. 2 (2019) 71-97.

\bibitem{CJY} T. Collins, A. Jacob, and S.T. Yau,
``(1,1)-forms with specified Lagrangian phase: A priori estimates and algebraic obstructions'', Camb. J. Math. 8 (2020), no. 2, 407-452.

\bibitem{DeG} E. De Giorgi,  ``Sulla differenziabilit\`a e l'analiticit\`a delle estremali degli integrali multipli regolari''. Mem. Accad. Sci. Torino. Cl. Sci. Fis. Mat. Nat. (3) 3 1957 25 - 43.

\bibitem{DP} J.P. Demailly and N. Pali, ``Degenerate complex Monge-Amp\`ere equations over compact K\"ahler manifolds",
Intern. J. Math. 21 (2010) no. 3, 357-405.

\bibitem{DGL} E. Di Nezza, V. Guedj, and C.H. Lu, ``Finite energy vs finite entropy", arXiv: 2006.07061.

\bibitem{DK1}S. Dinew and S. Kolodziej, ``Pluripotential estimates on compact Hermitian Manifolds'', Advances in geometric analysis, (2012), 69 - 86, Adv. Lect. Math. (21).

\bibitem{DK2} S. Dinew and S. Kolodziej, ``A priori estimates for complex Hessian equations", Anal. PDE 7 no 1 (2013) 227-244.

\bibitem{D} S.K. Donaldson, ``Two-forms on four-manifolds and elliptic equations'', Inspired by S. S. Chern (World Scientific, 2006).

\bibitem{EGZ} P. Eyssidieux, V. Guedj, and A. Zeriahi, 
``Singular K\"ahler-Einstein metrics", J. Amer. Math. Soc. 22 (2009) 607-639.

\bibitem{FGS}X. Fu, B. Guo and J. Song, ``Geometric estimates for complex Monge-Amp\`ere equations'', J. Reine Angew. Math. 765 (2020), 69 - 99.


\bibitem{FWW} J.X. Fu, Z.Z. Wang, and D. Wu,
``Form-type Calabi-Yau equations", Math. Res. Lett. 17 (2010) 887-903.

\bibitem{FY} J.X. Fu and S.T. Yau,  ``The theory of superstring with flux on non-K\"ahler manifolds and the complex Monge-Ampère equation''. J. Differential Geom. 78 (2008), no. 3, 369 - 428.

\bibitem{GGQ} M. George, B. Guan, and C. Qiu,
``Fully nonlinear elliptic equations on Hermitian manifolds for symmetric functions of partial Laplacians", arXiv: 2110.00490.

\bibitem{G} B. Guan, ``Second-order estimates and regularity for fully nonlinear elliptic equations on Riemannian manifolds''. Duke Math. J. 163 (2014), no. 8, 1491 - 1524.

\bibitem{GN} B. Guan and X.L. Nie, ``Second order estimates for fully non-linear equations with gradient terms on Hermitian manifolds", arXiv: 2108.03308.


\bibitem{GL} V. Guedj and H.C. Lu, ``Quasi-plurisubharmonic envelopes 3: Solving Monge-Amp\`ere equations on hermitian manifolds", arXiv:2107.01938.


\bibitem{GPT} B. Guo, D.H. Phong, and F. Tong, ``On $L^\infty$ estimates for complex Monge-Amp\`ere equations'',   arXiv:2106.02224

\bibitem{GPT21s}B. Guo, D.H. Phong, and F. Tong, ``A new gradient estimate for the complex Monge-Amp\`ere equation'', preprint, arXiv:2106.03308

\bibitem{GPT1} B. Guo, D.H. Phong, and F. Tong, ``Stability estimates for the complex Monge-Amp\`re and Hessian equations'', to appear in Calc. Var. Partial Differ. Equ.  arXiv:2106.03913

\bibitem{GPTW1} B. Guo, D.H. Phong, F. Tong, and C. Wang, ``On $L^\infty$ estimates for Monge-Amp\`ere and Hessian equations on nef classes'', to appear in Anal. PDE arXiv:2111.14186

\bibitem{GPTW2}B. Guo, D.H. Phong, F. Tong, and C. Wang,  ``On the modulus of continuity of solutions to complex Monge-Amp\`ere equations'',  arXiv:2112.02354

\bibitem{GPS} B. Guo, D.H. Phong, and J. Sturm, ``Green's functions and complex Monge-Amp\`ere equations", arXiv:2202.04715.

\bibitem{GP} B. Guo and D.H. Phong, ``On $L^\infty$ estimates for fully nonlinear partial differential equations on Hermitian manifolds'', preprint arXiv:2204.12549

\bibitem{GPa} B. Guo and D.H. Phong, ``Uniform entropy and energy bounds for fully non-linear equations'', preprint arXiv:2207.08983

\bibitem{GPSS} B. Guo, D.H. Phong, J. Song, and J. Sturm, ``Diameter estimates in K\"ahler geometry'', preprint arXiv:2209.09428

\bibitem{GS} B. Guo and J. Song,  ``Local noncollapsing for complex Monge-Amp\`ere equations'', to appear in {J. Reine Angew. Math.} arXiv:2201.02930

\bibitem{HL1} F. R. Harvey and H. B. Lawson,  ``Dirichlet duality and the nonlinear Dirichlet problem''. Comm. Pure Appl. Math. 62 (2009), no. 3, 396 - 443. 

\bibitem{HL2}  F. R. Harvey and H. B. Lawson,  ``Dirichlet duality and the nonlinear Dirichlet problem on Riemannian manifolds'', J. Differential Geom. 88 (2011), no. 3, 395 - 482.

\bibitem{HL3} F.R. Harvey and H.B. Lawson, ``Determinant majorization and the work of Guo-Phong-Tong and Abja-Olive", arXiv: 2207.01729.

\bibitem{H}L. H\"ormander, ``An introduction to complex analysis in several variables'',  Van Nostrand, Princeton, NJ, 1973


\bibitem{K} S. Kolodziej, ``The complex Monge-Amp\`ere equation", Acta Math. 180 (1998) 69-117.

\bibitem{Ku1} M. Kuranishi, ``Strongly pseudoconvex CR structures over small balls. I. An a priori estimate''.
Ann. of Math. (2) 115 (1982), no. 3, 451 - 500.

\bibitem{Ku2} M. Kuranishi, ``Strongly pseudoconvex CR structures over small balls. II. A regularity theorem''.
Ann. of Math. (2) 116 (1982), no. 1, 1 - 64.

\bibitem{Ku3}M. Kuranishi, ``Strongly pseudoconvex CR structures over small balls. III. An embedding theorem''. Ann. of Math. (2) 116 (1982), no. 2, 249 - 330.

\bibitem{LY} P. Li and S.T. Yau,
``On the parabolic kernel of the Schr\"odinger operator", Acta Math. 156 (1986), no. 3-4, 153 - 201

\bibitem{P} D.H. Phong, ``Geometric partial differential equations from unified string theories",  Proceedings of the International Consortium of Chinese Mathematicians 2018, 67–87, Int. Press, Boston, MA, 2020.

\bibitem{PPZ0} D.H. Phong, S. Picard, and X.W. Zhang,
``A second order estimate for general complex Hessian equations", Analysis \& PDE 9 no. 7 (2017) 1693-1709.

\bibitem{PPZ1} D.H. Phong, S. Picard, and X.W.  Zhang, ``Anomaly Flows",
Comm. Anal. Geom. 26 (2018) no. 4, 955-1008, arXiv: 1610.02739.


\bibitem{PPZ} D.H. Phong, S. Picard, and X.W. Zhang,
``Fu-Yau Hessian equations", J. Differential Geom. 118 (1) (2021) 147-187.

\bibitem{PSS}D.H. Phong, N. Sesum, and J. Sturm, ``Multiplier ideal sheaves and the K\"ahler-Ricci
flow'', Comm. Anal. Geom. 15 (2007), no. 3, 613 - 632.

\bibitem{PS09}D.H. Phong and J. Sturm, ``The Dirichlet problem for degenerate complex Monge-Amp\`ere equations'', Comm. Anal. Geom. 18 (2010), no. 1, 145 - 170.

\bibitem{Po} D. Popovici, ``Non-K\"ahler mirror symmetry of Iwasawa manifolds", arXiv: 1706.06449, to appear in Int. Math. Res. Notices.

\bibitem{SY} R. Schoen and S.T. Yau, ``Lectures on differential geometry''. Conference Proceedings and Lecture Notes in Geometry and Topology, I. International Press, Cambridge, MA, 1994. v+235 pp. 

\bibitem{Sl} L. Simon, ``Theorems on regularity and singularity of energy minimizing maps''. Lectures in Mathematics ETH Z\"urich. Birkh\"auser Verlag, Basel, 1996. viii+152 pp.

\bibitem{SWYY} I.M. Singer, B. Wong, S.T. Yau, and S.S.T. Yau, ``An estimate of the gap of the first two eigenvalues in the Schrödinger operator''. Ann. Scuola Norm. Sup. Pisa Cl. Sci. (4) 12 (1985), no. 2, 319 - 333.

\bibitem{ST} J. Song and G. Tian, ``Bounding scalar curvature for global solutions of the K\"ahler-Ricci flow''. Amer. J. Math. 138 (2016), no. 3, 683 - 695.

\bibitem{S} J. Spruck, ``Geometric aspects of the theory of fully non linear elliptic equations",  Global theory of minimal surfaces, 283 - 309, Clay Math. Proc., 2, Amer. Math. Soc., Providence, RI, 2005

\bibitem{Sz} G. Sz\'ekelyhidi  ``Fully non-linear elliptic equations on compact Hermitian manifolds''. J. Differential Geom. 109 (2018), no. 2, 337 - 378.

\bibitem{STW}G. Sz\'ekelyhidi, V. Tosatti and B. Weinkove ``Gauduchon metrics with prescribed volume form''. Acta Math. 219 (2017), no. 1, 181 - 211. 

\bibitem{T} G. Tian, ``On K\"ahler-Einstein metrics on certain K\"ahler manifolds with $C_1(M)>0$'',  Invent. Math. 89 (1987), no. 2, 225--246

\bibitem{ToWe} V. Tosatti, and B. Weinkove, ``The complex Monge-Amp\`ere equation on compact Hermitian manifolds''. J. Amer. Math. Soc. 23 (2010), no. 4, 1187 - 1195.

\bibitem{TW2} V. Tosatti, and B. Weinkove ``Hermitian metrics, $(n-1,n-1)$ forms and Monge-Amp\`ere equations''. J. Reine Angew. Math. 755 (2019), 67 - 101. 

\bibitem{TW1} V. Tosatti, and B. Weinkove ``The Monge-Amp\`ere equation for $(n-1)$-plurisubharmonic functions on a compact K\"ahler manifold'', J. Amer. Math. Soc. 30 (2017), no.2, 311-346.

\bibitem{TWY}V. Tosatti,  B. Weinkove and S.-T. Yau, ``Taming symplectic forms and the Calabi-Yau equation'', Proc. London Math. Soc. 97 (2008), no.2, 401- 424.

\bibitem{WWZ} J.X. Wang, X.J. Wang, and B. Zhou, ``A priori estimates for the complex Monge-Amp\`ere equation", arXiv:2003.06059.
\bibitem{W} X.J. Wang, ``Schauder estimates for elliptic and parabolic equations''. Chinese Ann. Math. Ser. B 27 (2006), no. 6, 637 - 642.

\bibitem{W} B. Weinkove, ``The Calabi-Yau equation on almost-K\"ahler four-manifolds''. J. Differential Geom. 76 (2007), no. 2, 317 - 349. 


\bibitem{Y} S.T. Yau, ``On the Ricci curvature of a compact K\"ahler manifold and the complex Monge-Amp\`ere equation. I'', Comm. Pure Appl. Math. 31 (1978) 339-411.



}

\end{thebibliography}
\end{document}